\documentclass[12pt]{article}

\usepackage{a4wide}
\usepackage{amssymb}
\usepackage{amsfonts}
\usepackage{amsmath}
\input xy
\xyoption{arrow}
\xyoption{matrix}

\date{}

\newtheorem{proposition}{Proposition}[section]
\newtheorem{theorem}[proposition]{Theorem}
\newtheorem{lemma}[proposition]{Lemma}

\newtheorem{corollary}[proposition]{Corollary}

\def\GK{{\rm  GK}\,}
\def\Kdim{{\rm K.dim }\,}

\def\Hom{{\rm Hom}}
\def\der{\partial }

\def\nFM0{{\nu }_{F,M_0}}
\def\nFN0{{\nu }_{F,N_0}}
\def\nGN0{{\nu }_{G,N_0}}

\def\N0{ {\bf N}_0 }

\def\t{\otimes}
\def\g{\gamma}

\def\ra{\rightarrow}

\def\lra{\leftrightarrow}
\def\Xpm{X^{\pm }}

\def\s{\sigma}

\def\l1{{\lambda}_1}

\def\a{\alpha}
\def\a0{ {\alpha }_0}
\def\a1{ {\alpha }_1}

\def\l{\lambda}

\def\nFGM0{{\nu }_{F,G,M_0}}

\def\nFN0{{\nu}_{F,N_0}}

\def\sm{{\sigma}^m}

\def\sm1{{\sigma}^{-1}}

\def\smtp1{{\sigma}^{-t+1}}

\def\S1{S^{-1}}

\def\Xpm1{X^{\pm 1}_1}

\def\sPM1{{\sigma }^{\pm 1}}
\def\sMP1{{\sigma }^{\mp 1 }}

\def\d{\delta}

\def\OO{{\cal O}}
\def\CA{{\cal A}}

\def\CD{{\cal D}}

\def\Ytm1{Y^{t-1}}
\def\Yim1{Y^{i-1}}

\def\CK{{\cal K}}

\def\i{{\bf i}}

\def\Aut{{\rm Aut}}

\def\Der{{\rm Der }}
\def\ad{{\rm ad }}
\def\dim{{\rm dim }}

\def\ker{ {\rm ker } }

\def\D{ \Delta }
\def\SL2Z{ {\rm SL}_2({\bf Z}) }

\def\CR{ {\cal R}}

\def\Gp1{ G^{1 , 1 } }
\def\P11{ P^{-1 , 1 } }
\def\Pp1{ P^{1 , 1 } }

\def\nCLsr{{}^\nu\kern-2pt {\cal L}^{\sigma , \rho  }}
\def\nP{{}^\nu \kern-2pt P}
\def\nL{{}^\nu\kern-2pt L}
\def\nLL{{}^\nu\kern-2pt \Lambda}
\def\nPsr{{}^\nu\kern-2pt P^{\sigma , \rho  }}
\def\nLsr{{}^\nu\kern-2pt L^{\sigma , \rho  }}
\def\nuCL{{}^\nu\kern-2pt  {\cal L}}
\def\nCLsr{{}^\nu\kern-2pt {\cal L}^{\sigma , \rho  }}
\def\nCL1m{{}^\nu\kern-2pt {\cal L}^{-1 , 1  }}
\def\x1nu{x^\frac{1}{\nu}}
\def\xm1nu{x^{-\frac{1}{\nu}}}

\def\CR{ {\cal R}}

\def\ra{\rightarrow }

\def\CI{{\cal I}}

\def\CC{ {\cal C}}

\def\CP{ {\cal P}}

\def\nAM0{{\nu }_{{\cal A},M_0}}
\def\nAN0{{\nu }_{{\cal A},N_0}}

\def\Kdim{ {\rm Kdim } }

\def\End{ {\rm End }}
\def\Der{ {\rm Der }}

\def\CR{ {\cal R }}
\def\CP{ {\cal P }}
\def\det{ {\rm det }}
\def\ad{ {\rm ad }}

\def\bp{\overline{p}}

\def\bx{\overline{x}}

\def\ga{\mathfrak{a}}

\def\gm{\mathfrak{m}}
\def\gp{\mathfrak{p}}

\def\gD{\mathfrak{D}}

\def\gder{\gd\ge\gr}
\def\bJ{\overline{J}}
\def\tb{\widetilde{b}}
\def\tB{\widetilde{B}}
\def\hj{\widehat{j}}
\def\hi{\widehat{i}}
\def\dij{\delta_{{\bf i}, {\bf j}}}

\def\vijnu{v_{{\bf i}, j_\nu}}

\def\vijmu{v_{{\bf i}, j_\mu}}

\def\vijk{v_{{\bf i}, j_k}}

\def\vijmus{v_{{\bf i}, j_\mu'}}

\def\Dijnu{\Delta_{{\bf i}, j_\nu}}

\def\Dijnus{\Delta_{{\bf i}, j_\nu'}}

\def\Disj{\Delta_{{\bf i}', {\bf j}}}

\def\Hjjs{H( {\bf j}, {\bf j}')}
\def\Hjsj{H( {\bf j}', {\bf j})}

\def\derij{\partial_{{\bf i}, {\bf j}}}
\def\GL{\rm GL}
\def\j{{\bf j}}
\def\bJisjs{\overline{J}({\bf i}', {\bf j}')}

\def\Spec{{\rm Spec}}

\def\Hom{{\rm Hom}}
\def\II{{\bf I}}
\def\JJ{{\bf J}}
\def\gder{\mathfrak{der}}
\begin{document}

\author{V. V. \  Bavula}

\title{Generators and defining relations for
the ring of differential operators on a smooth affine algebraic
variety}

\maketitle
\begin{abstract}
For the ring of differential operators on a smooth affine
algebraic variety $X$ over a field of characteristic zero a {\em
finite} set of algebra generators and a {\em finite} set of {\em
defining} relations are found {\em explicitly}. As a consequence,
a finite set of generators and a finite set of defining relations
are given for the module $\Der_K(\OO (X))$ of derivations on the
algebra $\OO (X)$ of regular functions  on the variety $X$. For
the variety $X$ which is not necessarily smooth, a set of {\em
natural} derivations ${\rm der}_K(\OO (X))$ of the algebra $\OO
(X)$  and a ring $\gD (\OO (X))$ of {\em natural} differential
operators on $\OO (X)$ are introduced. The algebra $\gD (\OO (X))$
is a  Noetherian algebra of Gelfand-Kirillov dimension $2\dim
(X)$. When $X$ is smooth then ${\rm der}_K(\OO (X))=\Der_K(\OO
(X))$ and $\gD (\OO (X))=\CD (\OO (X))$. A criterion of smoothness
of $X$ is given when $X$ is irreducible ($X$ is smooth iff $\gD
(\OO (X))$ is a simple algebra iff $\OO (X)$ is a simple $\gD (\OO
(X))$-module). The same results are true for regular algebras of
 essentially finite type. For a  singular irreducible affine
algebraic variety $X$, in general, the algebra of differential
operators $\CD (\OO (X))$ needs not be finitely generated nor
(left or right) Noetherian, it is proved that each term $\CD (\OO
(X))_i$ of the  order filtration  $\CD (\OO (X))=\cup_{i\geq 0}\CD
(\OO (X))_i$ is a  finitely generated left $\OO (X)$-module.

\end{abstract}

%%%%%%%%%%%%%%%%%% SECTION 1 %%%%%%%%%%%%%%%%%%%%%%%%

\section{Introduction}

The following notation will remain  fixed throughout the paper (if
it is not stated otherwise): $K$ is a field of characteristic zero
(not necessarily algebraically closed), module means a left
module, $P_n=K[x_1, \ldots , x_n]$ is a polynomial algebra over
$K$, $\der_1:=\frac{\der}{\der x_1}, \ldots ,
\der_n:=\frac{\der}{\der x_n}\in \Der_K(P_n)$, $I:=\sum_{i=1}^m
P_nf_i$ is a {\bf prime} but {\bf not} a maximal ideal of the
polynomial algebra $P_n$ with a set of generators $f_1, \ldots ,
f_m$, the algebra $A:=P_n/I$ which is a domain with the field of
fractions $Q:={\rm Frac}(A)$, the epimorphism $\pi :P_n\ra A$,
$p\mapsto \bp :=p+I$, to make notation simpler we sometime write
$x_i$ for $\overline{x}_i$ (if it does not lead to confusion), the
{\bf Jacobi} $m\times n$ matrices
 $J=(\frac{\der f_i}{\der x_j})\in M_{m,n}(P_n)$ and
 $\bJ =(\overline{\frac{\der f_i}{\der
x_j}})\in M_{m,n}(A)\subseteq M_{m,n}(Q)$,  $r:={\rm rk}_Q(\bJ )$
is the {\bf rank} of the Jacobi matrix $\bJ$ over the field $Q$,
$\ga_r$ is the {\bf Jacobian ideal} of the algebra $A$ which is
(by definition) generated by all the $r\times r$ minors of the
Jacobi matrix $\bJ$ ($A$ is {\em regular} iff $\ga_r=A$, it is the
{\bf Jacobian criterion of regularity}, \cite{Eisenbook}, 16.19),
$\Omega_A$ is the module of {\bf K\"{a}hler} differentials for the
algebra $A$. For $\i =(i_1, \ldots , i_r)$ such that $1\leq
i_1<\cdots <i_r\leq m$ and $\j =(j_1, \ldots , j_r)$ such that
$1\leq j_1<\cdots <j_r\leq n$, $\D
 (\i , \j )$ denotes the corresponding minor of the Jacobi matrix
$\bJ =(\bJ_{ij})$, that is $\det (\bJ_{i_\nu, j_\mu})$, $\nu , \mu
=1, \ldots, r$,  and the $\i$ (resp. $\j $) is called {\bf
non-singular} if $\D (\i , \j')\neq 0$ (resp. $\D (\i', \j )\neq
0$) for some $\j'$ (resp. $\i'$). We denote by $\II_r$ (resp.
$\JJ_r$) the set of all the non-singular $r$-tuples $\i$ (resp.
$\j $).

Since $r$ is the rank of the Jacobi matrix $\bJ$, it is easy to
show that $\D (\i , \j )\neq 0$ iff $\i\in \II_r$ and $\j\in
\JJ_r$ (Lemma \ref{IIr1}). Denote by $\JJ_{r+1}$ the set of all
$(r+1)$-tuples $\j =(j_1, \ldots , j_{r+1})$ such that $1\leq
j_1<\cdots <j_{r+1}\leq n$ and when deleting  some element, say
$j_\nu$, we have a non-singular $r$-tuple $(j_1, \ldots
,\widehat{j_\nu},\ldots , j_{r+1})\in \JJ_r$ where the hat over
 a symbol here and everywhere means that the symbol is omitted from
the list. The set $\JJ_{r+1}$ is called the {\bf critical set} and
any element of it is called a {\bf critical singular}
$(r+1)$-tuple. $\Der_K(A)$ is the $A$-module of $K$-derivations of
the algebra $A$. The action of a derivation $\d$ on an element $a$
is denoted by $\d (a)$.

 The next theorem gives a finite set of generators and a
finite set of  defining relations for the left $A$-module
$\Der_K(A)$ when $A$ is a  regular algebra.
\begin{theorem}\label{9bFeb05}%\marginpar{9bFeb05}
Let the algebra $A$ be a regular algebra. Then the left $A$-module
$\Der_K(A)$ is generated by the derivations $\der_{\i , \j }$, $\i
\in \II_r$, $\j \in \JJ_{r+1}$, where
\begin{eqnarray*}
 \derij  =  \der_{i_1, \ldots , i_r; j_1, \ldots , j_{r+1}}:= \det
 \begin{pmatrix}
  \overline{\frac{\der f_{i_1}}{\der x_{j_1}}} & \cdots &  \overline{\frac{\der f_{i_1}}{\der
  x_{j_{r+1}}}}\\
  \vdots & \vdots & \vdots \\
  \overline{\frac{\der f_{i_r}}{\der x_{j_1}}} & \cdots &  \overline{\frac{\der f_{i_r}}{\der
  x_{j_{r+1}}}}\\
  \der_{j_1}& \cdots & \der_{j_{r+1}}\\
\end{pmatrix}
\end{eqnarray*}

that satisfy the following defining relations (as a left
$A$-module): %\marginpar{Derel}
\begin{equation}\label{Derel}
\D (\i , \j )\der_{\i', \j'}=\sum_{l=1}^s(-1)^{r+1+\nu_l}\D (\i';
j_1', \ldots , \widehat{j_{\nu_l}'}, \ldots , j_{r+1}')\der_{\i;\j
, j_{\nu_l}'}
\end{equation}
for all $\i, \i'\in \II_r$, $\j=(j_1, \ldots , j_r)\in \JJ_r$, and
$\j'=(j_1', \ldots , j_{r+1}')\in \JJ_{r+1}$ where $\{ j_{\nu_1}',
\ldots , j_{\nu_s}'\}=\{ j_1', \ldots , j_{r+1}'\}\backslash \{
j_1, \ldots , j_r\}$.
\end{theorem}

The next result gives a finite set of generators and a finite set
of  defining relations for the $K$-algebra $\CD (A)$ of
differential operators on $A$.

\begin{theorem}\label{9Feb05}%\marginpar{9Feb05}
Let the algebra $A$ be a regular algebra. Then the ring of
differential operators $\CD (A)$ on $A$ is generated over $K$ by
the algebra $A$ and the derivations  $\der_{\i , \j }$, $\i \in
\II_r$, $\j \in \JJ_{r+1}$ that satisfy the defining relations
(\ref{Derel}) and %\marginpar{1Derel} 
\begin{equation}\label{1Derel}
\der_{\i , \j }\overline{x}_k=\overline{x}_k\der_{\i , \j
}+\der_{\i , \j }(\overline{x}_k), 
 \der_{\i ; \tilde{\j}, \tilde{j}_s}(\D)\der_{\i;\tilde{\j}, \tilde{j_t}}=\der_{\i ; \tilde{\j}, \tilde{j}_t }(\D)\der_{\i;\tilde{\j}, \tilde{j}_s},
\i \in \II_r,  \j \in
\JJ_{r+1}, \tilde{\j}\in
\JJ_r,  1\leq k\leq n,
\end{equation}
$s,t\in \{ r+1, \ldots , n\}$, $s\neq t$ where $\D = \D (\i , \tilde{\j})$, $\tilde{\j}=(\tilde{j}_1, \ldots , \tilde{j}_r)$ and $\{ \tilde{j}_{r+1}, \ldots , \tilde{j}_n\} = \{ 1, \ldots , n\}\backslash \{ \tilde{j}_1, \ldots , \tilde{j}_r\}$.
More formally, as an abstract $K$-algebra the ring of differential
operators $\CD (A)$ is generated by elements $x_1, \ldots , x_n$
and $d_{\i , \j }$, $\i \in \II_r$, $\j \in \JJ_{r+1}$ that
satisfy the following defining relations (where $s,t\in \{ r+1, \ldots , n\}$, $s\neq t$):

%\marginpar{RD1}
\begin{equation}\label{RD1}
f_1=\cdots =f_m=0, \;\; x_ix_j=x_jx_i, \; i,j=1, \ldots , n,
\end{equation}
%\marginpar{RD2}
\begin{equation}\label{RD2}
d_{\i , \j }x_k=x_kd_{\i , \j }+\der_{\i , \j }(x_k), \der_{\i ; \tilde{\j}, \tilde{j}_s}(\D)\der_{\i;\tilde{\j}, \tilde{j_t}}=\der_{\i ; \tilde{\j}, \tilde{j}_t }(\D)\der_{\i;\tilde{\j}, \tilde{j}_s},
\i \in \II_r,  \j \in
\JJ_{r+1}, \tilde{\j}\in
\JJ_r,  1\leq k\leq n,
\end{equation}
%\marginpar{RD3}
\begin{equation}\label{RD3}
\D (\i , \j )d_{\i', \j'}=\sum_{l=1}^s(-1)^{r+1+\nu_l}\D (\i';
j_1', \ldots , \widehat{j_{\nu_l}'}, \ldots , j_{r+1}')d_{\i;\j ,
j_{\nu_l}'}
\end{equation}
for all $\i, \i'\in \II_r$, $\j=(j_1, \ldots , j_r)\in \JJ_r$, and
$\j'=(j_1', \ldots , j_{r+1}')\in \JJ_{r+1}$ where $\{ j_{\nu_1}',
\ldots , j_{\nu_s}'\}=\{ j_1', \ldots , j_{r+1}'\}\backslash \{
j_1, \ldots , j_r\}$, and, for each $\i \in \II_r$, $\j'\in
\JJ_{r+1}$, and for a permutation $\s \in S_{r+1}$, the element
$d_{\i; j_{\s (1)}', \ldots j_{\s (r+1)'}}$ means $\epsilon (\s
)d_{\i , \j'}$ where $\epsilon (\s )$ is the parity of $\s $.
\end{theorem}
%{\it Remark}.  The element $\der_{\i , \j}(\bx_k)$ in
%(\ref{1Derel}) (and the element $\der_{\i , \j}(x_k)$ in
%(\ref{RD2}))  means $(-1)^{r+1+s}$ $  \D (\i ; j_1, \ldots ,
%\widehat{j_s}, \ldots , j_{r+1})$ if $k=j_s$ for some $s$ where
%$\j =(j_1, \ldots , j_{r+1})$, and zero otherwise.

%It is well-known fact that the ring of differential operators $\CD
%(A)$  on a {\em regular} algebra $A$ over the field $K$ of
%characteristic {\em zero} is generated by $A$ and $\Der_K(A)$.

 The algebra $A$ is the algebra of regular functions on the irreducible
 affine algebraic variety $X=\Spec (A)$, therefore we have the
 explicit algebra generators for the ring of differential
 operators $\CD (X)= \CD (A)$ on an arbitrary  smooth irreducible affine
 algebraic variety $X$. Any regular affine algebra $A'$ is a
 finite direct product of regular affine domains,
 $A'=\prod_{i=1}^s A_i$.  {\em  Since $\CD (A')\simeq \prod_{i=1}^s
 \CD (A_i)$, Theorem \ref{9Feb05} gives  algebra generators and defining relations for the ring of
 differential operators on arbitrary smooth affine algebraic
 variety. Since $\Der_K(A')\simeq \bigoplus_{i=1}^s\Der_K(A_i)$,
 Theorem \ref{9bFeb05} gives generators and defining relations for the left
 $A'$-module of derivations} $\Der_K(A')$.

Recall an important criterion of regularity for the algebra $A$
via properties of the {\bf derivation} algebra $\D (A)$, it is a
subalgebra of the ring of differential operators $\CD (A)$
generated by the algebra $A$ and the derivations $\Der_K(A)$.

\begin{theorem}\label{MR15.3.8}%\marginpar{MR15.3.8}
(Criterion of regularity via $\D (A)$, \cite{MR}, 15.3.8) The
following statements are equivalent.
\begin{enumerate}
\item $A$ is a regular algebra. \item $\D (A)$ is a simple
algebra.\item $A$ is a simple $\D (A)$-module.
\end{enumerate}
\end{theorem}

In the general case (when $A$ is not necessarily regular) the
$A$-module ${\rm der}_K(A):=\sum_{\i \in \II_r, \j\in
\JJ_{r+1}}A\derij$ is called the $A$-module of {\bf natural
derivations} of $A$.  The algebra of {\bf natural differential
operators} $\gD (A)$ is a subalgebra of $\CD (A)$ generated by $A$
and ${\rm der }_K(A)$. It is always a (left and right) Noetherian
algebra of Gelfand-Kirillov dimension $2\Kdim (A)$ (Lemma
\ref{KdA}). The left $A$-module ${\rm der }_K(A)$ and the algebra
$\gD (A)$ {\em do not} depend on the presentation of the algebra
$A$ as a factor algebra $P_n/I$ (Theorem \ref{7Mar05}).
 Clearly, $\gD (A)\subseteq \D (A)\subseteq
\CD (A)$. In general, the inclusions are strict ({\em eg}, for the
cusp $y^2=x^3$). If the algebra $A$ is {\em regular} then $\gD
(A)= \D (A)=\CD (A)$ (Theorem \ref{9Feb05}). The next result is a
similar criterion of regularity as Theorem \ref{MR15.3.8} but
using properties of the algebra $\gD (A)$. It is a practical
criterion as we know explicitly the generators of the algebra $\gD
(A)$.

\begin{theorem}\label{25Feb05}%\marginpar{25Feb05}
(Criterion of regularity via $\gD (A)$) The following statements
are equivalent.
\begin{enumerate}
\item $A$ is a regular algebra. \item $\gD (A)$ is a simple
algebra.\item $A$ is a simple $\gD (A)$-module.
\end{enumerate}
\end{theorem}

{\it Remark}. It is an interesting fact that the coefficients in
the defining relations (\ref{Derel}) and the coefficients of the
$\der_{j_\nu}$ when one expands the determinant $\derij$  are
generators for the Jacobian ideal $\ga_r$ that determines the {\em
singular locus} of the variety $X=\Spec (A)$.

If the algebra $A$ is regular then  the ring of differential
operators $\CD (A)$  is a finitely generated Noetherian algebra.
If $A$ is  not regular  then, in general, the algebra $\CD (A)$
need not be  a finitely generated algebra nor a left or right
Noetherian algebra, \cite{BGGDiffcone72}, the algebra $\CD (A)$
can be finitely generated and right Noetherian yet not left
Noetherian, \cite{SmStafDifopcurve}.  Though, a kind of finiteness
still holds for a {\em singular} algebra $A$.

\begin{theorem}\label{CDAifg}%\marginpar{CDAifg}
Let $\CD (A)=\cup_{i\geq 0}\CD (A)_i$ be the order filtration of
$\CD (A)$. Then, for each $i\geq 0$, $\CD (A)_i$ is a finitely
generated left $A$-module.
\end{theorem}

The paper is organized as follows. Section \ref{DerPnIS} contains
technical results that are used  throughout  the paper. The first
parts of Theorems \ref{9bFeb05} and \ref{9Feb05} are proved
(Theorem \ref{23July04}). It is shown that the module of natural
derivations ${\rm der}_K(A)$ of $A$ and the ring $\gD (A)$ of
natural differential operators on $A$ are {\em invariants} for the
algebra $A$, i.e. they do not depend on the choice of generators
for the ideal $I$ and on the choice of a presentation of the
algebra $A$ as the factor algebra $P_n/I$ (Theorem \ref{7Mar05}).
Similar results are true for the modules of `higher' natural
derivations ${\rm der}_K(A)_s$, $s=1, \ldots , r$ (Theorem
\ref{8Mar05}). Theorem \ref{25July04} describes explicitly the
$A$-module of derivations $\Der_K(A)$. It is proved that the {\em
tangent space} $T(X)_x:= \Der_K(A_\gm , k_\gm )$ at a {\em simple}
point $x=x_\gm$ of the variety $X=\Spec (A)$ is canonically
isomorphic to $\Der_K(A_\gm )/\gm \Der_K(A_\gm )$ (Corollary
\ref{c325July04}), and that for a singular point, in general, this
is not the case (even in the case of an isolated singularity, the
cusp, $y^2=x^3$ at $x=y=0$).

In Section \ref{Reldefops}, Theorems \ref{9bFeb05}, \ref{9Feb05}
and \ref{25Feb05} are proved. As a byproduct, a short elementary
direct proof is given of the fact that the algebra $\CD (A)$ of
differential operators on a regular algebra $A$ is a {\em simple
Noetherian}  algebra generated by $A$ and $\Der_K(A)$ (Theorem
\ref{3Feb05}). For an arbitrary (not necessarily regular) algebra
$A$, the Gelfand-Kirillov dimension $\GK (\gD (A))=2\Kdim (A)$
(the Krull dimension) and the (left and right) Krull dimension
$\Kdim (A)\leq \Kdim (\gD (A))\leq 2\Kdim (A)$ (Lemma \ref{KdA}).

In Section \ref{essgrel}, the same results as Theorems
\ref{9bFeb05}, \ref{9Feb05}, and \ref{25Feb05} are proved but for
regular algebras of essentially finite type (Theorems
\ref{ess9bFeb05}, \ref{ess9Feb05}, and \ref{ess25Feb05}).

In Section \ref{opsing}, Theorem \ref{CDAifg} is proved
(Proposition \ref{sinDAfg}).

%%%%%%%%%%%%%%%%%% SECTION 2 %%%%%%%%%%%%%%%%%%%%%%%%

\section{Generators for the ring of differential operators
 on a smooth irreducible affine algebraic
variety}\label{DerPnIS}%\marginpar{DerPnIS}

In this section, the $A$-module of natural derivations ${\rm
der}_K(A)$, the $A$-module of exceptional derivations and the ring
of natural differential operators $\gD (A)$ are introduced. Their
definitions are given for the presentation
 of the algebra $A$ as the factor algebra $P_n/I$ and for a set of
 generators of the ideal $I$. It will be proved that, in fact, ${\rm
der}_K(A)$ and $\gD (A)$ are invariants for the algebra $A$ (they
do not depend on the presentation and on the choice of generators,
Theorem \ref{7Mar05}). Similar results are true for {\em higher}
natural derivations ${\rm der}_K(A)_s$ that will be introduced at
the end of this section.

Let $B$ be a commutative $K$-algebra. The ring of ($K$-linear)
{\bf differential operators} $\CD (B)$ on $B$ is defined as a
union of $B$-modules  $\CD (B)=\cup_{i=0}^\infty \,\CD_i (B)$
where $\CD_0 (B)=\End_R(B)\simeq B$, ($(x\mapsto bx)\lra b$),
$$ \CD_i (B)=\{ u\in \End_K(B):\, [r,u]:=ru-ur\in \CD_{i-1} (B)\; {\rm for \; each \; }\; r\in B\}.$$
 The set of $B$-modules $\{ \CD_i (B)\}$ is the {\bf order filtration} for
the algebra $\CD (B)$:
$$\CD_0(B)\subseteq   \CD_1 (B)\subseteq \cdots \subseteq
\CD_i (B)\subseteq \cdots\;\; {\rm and}\;\; \CD_i (B)\CD_j
(B)\subseteq \CD_{i+j} (B), \;\; i,j\geq 0.$$

The subalgebra $\D (B)$ of $\CD (B)$ generated by $B\equiv
\End_R(B)$ and the set ${\rm Der}_K (B)$ of all $K$-derivations of
$B$ is called the {\bf derivation ring} of $B$.

Suppose that $B$ is a  regular affine  domain of Krull dimension
$n<\infty $. In geometric terms, $B$ is the coordinate ring $\OO
(X)$ of a smooth irreducible  affine algebraic variety $X$ of
dimension $n$. Then
\begin{itemize}
\item ${\rm Der}_K (B)$ {\em is a finitely generated projective}
$B$-{\em module of rank} $n$, \item  $\CD (B)=\Delta (B) $, \item
$\CD (B)$ {\em is a simple (left and right) Noetherian domain of
Gelfand-Kirillov dimension}  $\GK \, \CD (B)=2n$ ($n=\GK (B)=\Kdim
(B))$.
\end{itemize}

For the proofs of the statements above the reader is referred to
\cite{MR}, Chapter 15.
 So, the domain $\CD (B)$ is a simple finitely generated infinite dimensional Noetherian algebra
(\cite{MR}, Chapter 15).

 The {\em gradient}  map $\nabla =(\der_1,
\ldots , \der_n):P_n\ra P_n^n$, $p\mapsto (\frac{\der p}{\der
x_1}, \ldots , \frac{\der p}{\der x_n})$ satisfies the following
property: $\nabla (pq)=\nabla (p) q+p\nabla (q)$ for all $p,q\in
P_n$. It follows that $\CR =\CR_I:= (\nabla I+IP_n^n)/IP_n^n$ is
the $A$-submodule of $A^n$ generated by the rows of the Jacobi
matrix $\bJ$. The $A$-module $\CR_I$ is the {\em module of
relations} for the module of {\bf K\"{a}hler differentials}
$\Omega_A$: $0\ra \CR_I\stackrel{i}\rightarrow A^n\ra \Omega_A\ra
0$ is the short exact sequence of $A$-modules where $i$ is the
natural inclusion. The rank $r={\rm rk}_Q (\overline{\frac{\der
f_i}{\der x_j}})$ of the Jacobi matrix over the field $Q$ is equal
to $\dim_Q(Q\CR_I)$. $A^n=\bigoplus_{i=1}^nAe_i\subseteq
Q^n:=\bigoplus_{i=1}^nQe_i$ where $e_{i_1}=(1, 0, \ldots , 0),
\ldots , e_n=(0,\ldots , 0, 1)$. Inclusions of  $A$-modules $\CR
\subseteq A^n\subseteq Q^n$ induce $A$-module homomorphisms of the
exterior algebras over $A$: $\bigwedge^*\CR \ra
\bigwedge^*A^n\subseteq \bigwedge^*Q^n$. For each $k=1, \ldots ,
n$ and $1\leq i_1<\cdots <i_k\leq n$,  let $p_{i_1, \ldots ,
i_k}:\bigwedge^k\CR \ra A$ be the composition of two maps:
$\bigwedge^k \CR \ra \bigwedge^kA^n\ra \bigwedge^kQ^n$ and the
natural projection map $\bigwedge^kQ^n\ra Qe_{i_1}\wedge\cdots
\wedge e_{i_k}\simeq Q$. Then $\ga_k:= \sum_{1\leq i_1<\cdots
<i_k\leq n} p_{i_1, \ldots , i_k} (\bigwedge^k\CR)$ is the ideal
of $A$. {\em So, the ideal $\ga_k$ is generated by all the
$k\times k$ minors of the Jacobi matrix $(\overline{\frac{\der
f_i}{\der x_j}})$ and does not depend on the choice of the
generators of the ideal $I$}. In particular $\ga_{r+1}=\cdots =
\ga_n=0$, $\ga_1\supseteq \ga_2\supseteq \cdots \supseteq
\ga_r\neq 0$, $\ga_i\ga_j\supseteq \ga_{i+j}$ for all $i,j$. In
particular, $\ga_1^i\supseteq \ga_i$ for all $i$. The ideals
$\ga_i$ are called the {\bf Jacobian  ideals} for the algebra $A$.
It is well known that that Jacobian ideals $\ga_r, \ga_{r-1},
\ldots , \ga_1$ are {\em invariant} of the algebra $A=P_n/I$ in
the sense that if $ \s :A\simeq P_{n'}/I'$ is another presentation
for the algebra $A$ and $\ga_{r'}', \ga_{r'-1}', \ldots , \ga_1'$
are the corresponding Jacobian ideals then $\s (\ga_r)= \ga_{r'}',
\s (\ga_{r-1})= \ga_{r'-1}', \ldots $ (\cite{Eisenbook}, p. 405).
This is true because the Jacobian ideals are the {\bf Fitting
ideals} for the module of K\"{a}hler differentials $\Omega_A$:
$0\ra \CR_I\stackrel{i}\rightarrow A^n\ra \Omega_A\ra 0$, it is
well known that the Fitting ideals are invariants of a finitely
presented  module (\cite{Fitting1936}, see also \cite{Eisenbook},
 20.4).

% The Fitting ideals are invariant under the action of the group
% $\Aut_K(A)$ of $K$-algebra automorphisms of $A$: an automorphism
% $\s \in \Aut_K(A)$ determines new algebra generators for $A$,
% $\bx_1':=\s (\bx_1), \ldots , \bx_n':= \s (\bx_n)$ then
% $$ \bJ:= \frac{\der (\overline{f}_1, \ldots , \overline{f}_m)}{\der (\bx_1, \ldots ,
% \bx_n)}=\frac{\der (\overline{f}_1, \ldots , \overline{f}_m)}{\der
% (\bx_1', \ldots , \bx_n')}\frac{\der (\bx_1', \ldots ,
% \bx_n')}{\der (\bx_1, \ldots , \bx_n)}=\bJ'\frac{\der \bx'}{\der
% \bx}\;\; {\rm and}\;\; \bJ'=\bJ \frac{\der \bx}{\der \bx'},$$ and
% the result follows.

 Given an $n\times
n$ matrix $B=(B_{ij})$ with entries from an arbitrary field,
$1\leq i_1<\cdots <i_r\leq n $ and $1\leq j_1<\cdots <j_r\leq n $.
The determinant of the $r\times r$ matrix $(b_{ij})$ where
$i=i_1,\ldots , i_r$ and $j=j_1\ldots , j_r$ is called the {\em
minor} of the matrix $B$ denoted $\Delta (i_1,\ldots , i_r;
j_1\ldots , j_r)$. The $n\times n$ matrix $\tB=(\tb_{ij})$,
$\tb_{ij}:=(-1)^{i+j}\Delta (1, \ldots , \hj , \ldots , n;1,
\ldots , \hi , \ldots , n)$ is called the {\em adjoint matrix} to
$B$ or the {\em matrix of cofactors} of $B$ where the hat over a
symbol means that it is missed on the list. The adjoint matrix has
the property that $\tB B=B\tB =\det (B)E$ where $E$ is the
identity matrix. In particular, $\det (B) \d_{i,j}
=\sum_{k=1}^n\tb_{ik}b_{kj}=\sum_{k=1}^nb_{ik}\tb_{kj}$ for all
$i,j=1, \ldots , n$, where $\d_{i,j}$ is the Kronecker delta.

For each $r$-tuple $\i=(i_1, \ldots , i_r)$  where $i_1, \ldots ,
i_r\in \{ 1, \ldots , m\} $ and $(r+1)$-tuple $\j = (j_1,\ldots ,
j_{r+1})$ where $j_1,\ldots , j_{r+1}\in \{ 1, \ldots , n\}$,
consider the $K$-derivation of the polynomial algebra $P_n$ given
by the rule: %\marginpar{dbibj}
\begin{eqnarray}\label{dbibj}
 \dij  = \d_{i_1, \ldots , i_r; j_1, \ldots , j_{r+1}}:= \det
 \begin{pmatrix}
  \frac{\der f_{i_1}}{\der x_{j_1}} & \cdots &  \frac{\der f_{i_1}}{\der
  x_{j_{r+1}}}\\
  \vdots & \vdots & \vdots \\
  \frac{\der f_{i_r}}{\der x_{j_1}} & \cdots &  \frac{\der f_{i_r}}{\der
  x_{j_{r+1}}}\\
  \der_{j_1}& \cdots & \der_{j_{r+1}}\\
\end{pmatrix}
\end{eqnarray}
where the determinant means
\begin{eqnarray*}
\dij &=&\sum_{k=1}^{r+1}(-1)^{r+1+k}\D (i_1, \ldots , i_r; j_1,
\ldots ,
\widehat{j_k}, \ldots , j_{r+1})\der_{j_k}\\
  & =&  \D (i_1, \ldots , i_r; j_1, \ldots , j_{r})\der_{j_{r+1}}-
  \sum_{k=1}^r\D (i_1, \ldots , i_r; j_1, \ldots
, j_{k-1},j_{r+1}, j_{k+1}, \ldots , j_r)\der_{j_k} \\
\end{eqnarray*}
where $\der_i:=\frac{\der}{\der x_i}\in \Der_K(P_n)$. The
definition of $\d_{\i , \j }$ makes sense as the Krull dimension
of the algebra $A$ is $\Kdim (A)=n-r>0$, and so   $r<n$ and $r\leq
m$. If the algebra $A$ were a field then $\d_{\i , \j }$ would
make no sense.  Let $S_r$ be the symmetric group of order $r$ and
$\epsilon (\s )$ be the sign of a permutation $\s \in S_r$. For
$\s \in S_r$ and $\tau \in S_{r+1}$,
 $\d_{i_{\s (1)}, \ldots , i_{\s (r)}; j_{\tau (1)}, \ldots ,
j_{\tau (r+1)}}=\epsilon (\s ) \epsilon (\tau ) \d_{i_1 , \ldots ,
i_r; j_1, \ldots , j_{r+1}}$. Since (recall that $r$ is the rank
of the Jacobi matrix $\bJ$) %\marginpar{dbifk0}
\begin{equation}\label{dbifk0}
\dij (f_k)\equiv 0 \; \mod I, \;\; k=1, \ldots , m,
\end{equation}
we have $\dij (I)\subseteq I$ for all $\i $ and $\j $. So, each
derivation $\dij$ of the polynomial algebra $P_n$ induces the
$K$-derivation $\derij$ of the algebra $A:=P_n/I$ which is given
by the rule
 $\derij (p+I):= \dij (p)+I$. One can write %\marginpar{1dbibj}
\begin{eqnarray}\label{1dbibj}
 \derij  = \der_{i_1, \ldots , i_r; j_1, \ldots , j_{r+1}}:= \det
 \begin{pmatrix}
  \overline{\frac{\der f_{i_1}}{\der x_{j_1}}} & \cdots &  \overline{\frac{\der f_{i_1}}{\der
  x_{j_{r+1}}}}\\
  \vdots & \vdots & \vdots \\
  \overline{\frac{\der f_{i_r}}{\der x_{j_1}}} & \cdots &  \overline{\frac{\der f_{i_r}}{\der
  x_{j_{r+1}}}}\\
  \der_{j_1}& \cdots & \der_{j_{r+1}}\\
\end{pmatrix}.
\end{eqnarray}

{\it Example}. Let $I=(x_1, \ldots , x_r)$ be an ideal of the
polynomial algebra $P_n=K[x_1, \ldots , x_n]$, $r<n$. Then
$A=P_n/I\simeq K[x_{r+1}, \ldots , x_n]$ is a polynomial algebra
and $\der_{\i ; \i , k}=\frac{\der }{\der x_k}$, $k=r+1, \ldots ,
n$ where $\i :=(1, \ldots , r)$.

The next lemma is a result of linear algebra and an easy (but
rather technical) consequence of the fact that $r$ is the rank of
the (Jacobi) matrix.
\begin{lemma}\label{IIr1}%\marginpar{IIr1}
$\i \in \II_r$ and $\j\in \JJ_r$ $\Leftrightarrow$ $\D (\i , \j
)\neq 0$.
\end{lemma}
A proof will be given after Lemma \ref{nscoo}. We do this in order
to avoid unnecessary technicalities at this stage of the paper and
to make clear what is really important in finding the generators
for the ring of differential operators $\CD (A)$.

It follows from the definition of the derivations $\derij$ and
Lemma \ref{IIr1} that $\derij \neq 0$ iff $\i \in \II_r$ and $\j
\in \JJ_{r+1}$.

{\it Definition}. For the algebra $A=P_n/I$ and a given set $f_1,
\ldots , f_m$ of generators for the ideal $I$, we denote by ${\rm
der}_K(A)$ the $A$-submodule of $\Der_K(A)$ generated by all the
derivations $\der_{\i , \j }$ (see (\ref{1dbibj})), then ${\rm
der}_K(A)=\sum_{\i\in \II_r, \j\in\JJ_{r+1}}A\derij$ (by Lemma
\ref{IIr1}). We call ${\rm der}_K(A)$ the set of {\bf natural
derivations} of $A$, and an element of ${\rm der}_K(A)$ is called
a {\bf natural derivation} of $A$. A derivation of $A$ which is
not natural is called an {\bf exceptional derivation}, the left
$A$-module $\Der_K(A)/{\rm der}_K(A)$ is called the {\bf module of
exceptional derivations}. The algebra of {\bf natural differential
operators} $\gD (A)$ is the subalgebra of $\CD (A)$ generated by
$A$ and ${\rm der}_K(A)$.

{\it Example}. For the cusp, $A=K[x,y]/(x^3-y^2)$, we have
$\Der_K(A)=A\d +A\der $ and ${\rm der}_K(A)=A\d $ where $\d :=\det
 \begin{pmatrix}
  3x^2 & -2y\\
  \der_x&  \der_y\\
\end{pmatrix}=2y\der_x+3x^2\der_y$ and
$\der:=xy^{-1}\d=2x\der_x+3y\der_y$ (the Euler derivation). So,
the Euler derivation $\der $ is an exceptional derivation.

If the algebra $A$ is regular then ${\rm der}_K(A)=\Der_K(A)$
(Theorems \ref{9bFeb05} and  \ref{23July04}).

{\it Question}. {\em Does } ${\rm der}_K(A)=\Der_K(A)$ {\em imply
that $A$ is regular?}

The group $\Aut_K(A)$ of $K$-algebra automorphisms of the algebra
$A$ acts in a natural way (by changing generators) on the
derivations $\Der_K(A)$, $\s (\d ):=\s \circ \d \circ \s^{-1}$ for
$\s \in \Aut_K(A)$ and $\d \in \Der_K(A)$.

The derivations $\der_{\i , \j }$  depend on the choice of
generators of the ideal $I$. There are two natural questions: (1)
{\em Does ${\rm der}_K(A)$ depend on the choice of generators for
the ideal} $I$? (2) {\em Does ${\rm der}_K(A)$ depend on the
presentation $P_n/I$}? (i.e. given a $K$-algebra isomorphism $\s :
P_n/I\ra P_{n'}/I'$, is $\s {\rm der}_K(P_n/I)\s^{-1}={\rm
der}_K(P_{n'}/I')$?) The next theorem gives affirmative answers to
both questions, and so the module of natural derivations ${\rm
der}_K(A)$ is an {\em invariant} of the algebra $A$.

\begin{theorem}\label{7Mar05}%\marginpar{7Mar05}
\begin{enumerate}
\item The $A$-module ${\rm der}_K(A)$ does not depend on the
choice of generators for the ideal $I$. \item The set of natural
derivations   ${\rm der}_K(A)$ does not depend on the presentation
of the algebra $A$ as the factor algebra $P_n/I$. \item The ring
$\gD (A)$ of natural differential operators on $A$ does not depend
either on the choice of generators for the ideal $I$ or on the
presentation of the algebra $A$ as the factor algebra $P_n/I$.
\end{enumerate}
\end{theorem}

{\it Proof}.  $1$.  Given another set, say $g_1, \ldots , g_l$, of
generators for the ideal $I$. Let ${\rm der}_K(A)'$ be the
corresponding submodule of $\Der_K(A)$ of $g$-natural derivations.
Then $f_i=\sum_{k=1}^la_{ik}g_k$ for some $a_{ik}\in P_n$, and so
$\nabla f_i=\sum_{k=1}^la_{ik}\nabla g_{k}$(mod $I$). For $\j
=(j_1, \ldots , j_{r+1})$ and $f_i$, let $\nabla_\j
f_i:=(\overline{\frac{\der f_i}{\der x_{j_1}}}, \ldots ,
\overline{\frac{\der f_i}{\der
x_{j_{r+1}}}})=\sum_{i'=1}^l\overline{a}_{ii'}\nabla_\j g_{i'}$
and let $\nabla_\j := (\der_{j_1}, \ldots , \der_{j_{r+1}})$. The
derivations $\der_{\i , \j }$ can be written as the exterior
product (the determinant) $\nabla_\j f_{i_1}\wedge \cdots  \wedge
\nabla_\j f_{i_r}\wedge \nabla_\j$, where $\i =(i_1, \ldots ,
i_r)$. Substituting the expression for each $\nabla f_{i_\nu}$ in
the product, we see that %\marginpar{dijdijp}
\begin{equation}\label{dijdijp}
\der_{\i , \j}=\sum_{\i'}\det (a_{\i, \i'})\der_{\i', \j}'
\end{equation}
where $\i'=(i_1', \ldots , i_r')$ runs through all the $r$-tuples
such that $1\leq i_1'<\cdots <i_r'\leq l$, $a_{\i, \i'}$ is the
$r\times r$ matrix with $(\nu, \mu )$-entry equal to
$\overline{a_{i_\nu , i_{\mu}'}}$, and $\der_{\i', \j}'$ is  the
derivation (\ref{1dbibj}) for the second choice of generators. It
follows that ${\rm der}_K(A)\subseteq {\rm der}_K(A)'$. By
symmetry, the opposite inclusion is also true which proves that
${\rm der}_K(A)={\rm der}_K(A)'$.

$2$. Given another presentation $P_{n'}/I'$ for the algebra $A$
where $P_{n'}:=K[x_1', \ldots , x_{n'}']$. Fix a commutative
diagram
$$\xymatrix{P_{n+n'}\ar[d]&=& P_n\t P_{n'}\ar[d]\\
P_n\ar[d] & &  P_{n'}\ar[d]\\
P_n/I&\simeq & P_{n'}/I'}$$
 The idea of the proof is to use the
first statement. The LHS of the diagram gives a natural  algebra
isomorphism $ P_{n+n'}/I''\simeq P_n/I$ where $I'':= (f_1, \ldots
, f_m, x_1'-g_1, \ldots , x_{n'}'-g_{n'})$ for some polynomials
$g_1, \ldots , g_{n'}\in P_n$.  Set $x_{n+1}:= x_1', \ldots ,
x_{n+n'}:=x_{n'}'$. Then $P_{n+n'}=K[x_1, \ldots , x_{n+n'}]$ and,
for the presentation $A\simeq P_{n+n'}/I''$,  the corresponding
Jacobi matrix $\bJ''\in M_{n+n'}(A)$ has the form $\begin{pmatrix}
 \bJ & 0\\
  * & E \\
\end{pmatrix}$ where $E$ is the identity $n'\times n'$ matrix. The
rank of this matrix is $r+n'$. It follows at once from the Laplace
identity for determinant and from the fact that derivation is
determined by its action on any algebra generating set that
elements $\derij'':=\der_{\i, n+1, \ldots, n+n'; \j , n+1, \ldots,
n+n'}$ where $\i \in \II_r$ and $\j \in \JJ_{r+1}$ are  generators
for the presentation $P_{n+n'}/I''$. The actions of the
derivations $\derij''$ and $\derij$ on the generators $\bx_1,
\ldots , \bx_n$ of the algebra $A=P_n/I=P_{n+n'}/I''$ coincide,
therefore $\derij''=\derij$, and so ${\rm der}_K(A)={\rm
der}_K(A)''$ (it is in fact a natural isomorphism). By symmetry,
we have ${\rm der}_K(A)'\simeq {\rm der}_K(A)''$, and the result
follows, where ${\rm der}_K(A)'$ and ${\rm der}_K(A)''$ are the
natural derivations for the presentations $P_{n'}/I'$ and
$P_{n+n'}/I''$ respectively.

$3$. Evident.  $\Box $

\begin{corollary}\label{2Feb05}%\marginpar{2Feb05}
\begin{enumerate}
 \item The set of natural
derivations ${\rm der}_K(A)$ (as a subset of $\Der_K(A)$)  is
invariant under the action of the group $\Aut_K(A)$ (i.e. ${\rm
der}_K(A)$ does not depend on the choice of generators for the
algebra $A$ induced by an automorphism). \item $\Der_K(A)/{\rm
der}_K(A)$ is an $\Aut_K(A)$-module.
\end{enumerate}
\end{corollary}

Clearly, ${\rm der}_K(A)\subseteq \Der_K(A, \ga_r)=\{\d \in
\Hom_K(A, \ga_r)\, | \, \d (ab)=\d (a)b+a\d (b)$ for all $a,b\in
A\}$.

{\it Question}. {\em  Which property (if any, i.e. ${\rm
der}_K(A)=\Der_K(A, \ga_r)$) does distinguish the set of natural
derivations ${\rm der}_K(A)$ in the set $\Der_K(A, \ga_r)$ of all
derivations from $A$ to the Jacobian ideal $\ga_r$?}

Recall that $\Der_K(P_n)=\bigoplus_{i=1}^nP_n\der_i$. For any
derivation $\d' \in \Der_K(A)$ one can find a derivation $\d
\in\Der_K(P_n)$ such that $\d' \pi =\pi \d$, equivalently $\d
(I)\subseteq I$ and  $\d' (\bp )=\overline{\d (p)}$ for all $p\in
P_n$. Any derivation $\d $ of the polynomial algebra $P_n$ such
that $\d (I)\subseteq I$ induces
 the derivation $\d'$ of the algebra $A$ by the rule above. So,
%\marginpar{DerPnI}
\begin{equation}\label{DerPnI}
\Der_K(A)=\Der_K(P_n/I)\simeq \frac{\{ \d \in \Der_K(P_n): \; \d
(I)\subseteq I\}}{I\Der_K(P_n)}\simeq \{ \sum_{i=1}^na_i\der_i\in
\bigoplus_{i=1}^n A\der_i: \; \bJ a=0\},
\end{equation}
where $\bJ=(\overline{\frac{\der f_i}{\der x_j}} )$ is the Jacobi
matrix and $a=(a_1 , \ldots , a_n)^T$ is a  column vector.

\begin{proposition}\label{arDerA}%\marginpar{arDerA}
\begin{enumerate}
\item  $\D (i_1, \ldots , i_r; j_1, \ldots ,
j_r)\Der_K(A)\subseteq \sum_k A \der_{i_1, \ldots , i_r; j_1,
\ldots , j_r, k}$ for all $1\leq i_1< \cdots < i_r\leq m$ and $
1\leq j_1< \cdots < j_r\leq n$ where $k$ runs through the set $\{
1, \ldots , n\} \backslash \{ j_1, \ldots , j_r\}$. If $\D (i_1,
\ldots , i_r; j_1, \ldots , j_r)\neq 0$ then the sum above is the
direct sum. \item $\ga_r\Der_K(A)\subseteq {\rm der}_K(A)\subseteq
(\sum_{i=1}^n\ga_r\der_i)\cap \Der_K(A)$.
\end{enumerate}
\end{proposition}

{\it Proof}. $1$. It suffices to prove the first statement in the
case when $i_1=1, \ldots , i_r=r$ and $j_1=1, \ldots , j_r=r$. By
(\ref{DerPnI}),  $\Der_K(A)=\{ \sum_{i=1}^na_i\der_i\, | \, \bJ
a=0\}$ where $a=(a_1, \ldots , a_n)^T\in A^n$ is a  column vector.
So, $\bJ a=0$ is the system of $m$ linear equations with unknowns
$\{ a_i\}$. The first $r$ equations of the system can be written
in the matrix form as $Bb=-Cc$ where $B=(\overline{\frac{\der
f_i}{\der x_j}})_{i,j=1,\ldots , r}$, $C=(\overline{\frac{\der
f_i}{\der x_j}})^{i=1,\ldots , r}_{j=r+1,\ldots , m}$, $b=(a_1,
\ldots , a_r)^T$, and $c=(a_{r+1}, \ldots , a_n)^T$. Multiplying
the equation $Bb=-Cc$
 by the adjoint matrix $\tB=(\tb_{ij})$ of $B$ one has $\det (B)b=\tB
Bb=-\tB Cc$.
 In the matrix form, the derivation $\d =\sum_{i=1}^n a_i\der_i\in
 \Der_K(A)$  can be written as $\d =\nabla_1b+\nabla_2c$ where
 $\nabla_1:=(\der_1, \ldots , \der_r)$ and  $\nabla_2:=(\der_{r+1}, \ldots ,
 \der_n)$ and the products are formal, that is no action of the
 derivations $\der_i$ is assumed here and till the end of the
 proof. This is done only to simplify notation in the rest of the proof that follows,
 but if the
 reader feels uncomfortable with this assumption then apply the
 transpose to the equalities below  and then the coefficients are on the left
 from derivations and no assumption is needed. Let $\D =\D(1, \ldots , r; 1, \ldots , r)$ and
 $\i =(1, \ldots , r)$.  Then $\D =\det (B)$ and
 \begin{eqnarray*}
  \D\d &=& \det (B)\d
 =\nabla_1\det (B)b+\nabla_2\det (B)c= (\nabla_2 \det (B)-\nabla_1\tB
 C)c\\
 &=&\sum_{i=r+1}^n a_i(\det (B)\der_i-\sum_{j,k=1}^r\tb_{jk}
 \overline{\frac{\der f_k}{\der x_i}} \der_j)\\
 &=& \sum_{i=r+1}^na_i
  (\D \der_i-\sum_{j=1}^r
  \D(\i ; 1, \ldots ,j-1, i, j+1, \ldots ,  r)\der_j)\\
  &=& \sum_{i=r+1}^na_i
  ((-1)^{r+1+r+1}\D \der_i+\sum_{j=1}^r
 (-1)^{r+1+j} \D(\i ; 1, \ldots ,\hj , \ldots ,  r, i)\der_j)\\
 &=&\sum_{i=r+1}^na_i \der_{1, \ldots , r; 1, \ldots ,  r,
 i},
\end{eqnarray*}
by (\ref{1dbibj}), as required. If $\D \neq 0$ then, by
(\ref{1dbibj}), the sum in the first statement is direct (an
alternative argument: the rank of the $A$-module $\Der_K(A)$ is
$\Kdim (A)=n-r$, the sum above contains $n-r$ summands, so it must
be a direct sum).

$2$. Since the set of all the nonzero minors $\D (i_1, \ldots ,
i_r; j_1,\ldots , j_r)$ is a generating set for the Jacobian ideal
$\ga_r$, the second statement follows from the first.  $\Box $

The following theorem proves the first parts of Theorems
\ref{9bFeb05} and \ref{9Feb05}.
\begin{theorem}\label{23July04}%\marginpar{23July04}
Suppose that the algebra $A$ is a regular algebra. Then
\begin{enumerate}
\item  $\Der_K(A)={\rm der}_K(A)$. \item The algebra of
differential operators $\CD (A)$ is generated by the algebra $A$
and the derivations $\derij$, $\i \in \II_r$, $\j \in \JJ_{r+1}$.
\end{enumerate}
\end{theorem}

{\it Proof}. $1$. The algebra $A$ is regular iff the Jacobian
ideal $\ga_r$ is equal to $A$. So, statement 1 follows immediately
from Proposition \ref{arDerA}.(2): $\Der_K(A)\subseteq {\rm
der}_K(A)\subseteq \Der_K(A)$.

$2$. It is a well-known fact that for a regular  affine
commutative algebra $A$, the algebra of ($K$-linear) differential
operators $\CD (A)$ is generated by $A$ and $\Der_K(A)$ (see
\cite{MR}, 15.5.6, and also Theorem \ref{3Feb05}). $\Box $

As a corollary of Proposition \ref{arDerA} we have a short direct
proof of the following well-known result \cite{MR}, 15.2.11.

\begin{corollary}\label{carDerA}%\marginpar{carDerA}
Suppose that the algebra $A=P_n/I$ is a regular algebra. Then
$\Der_K(A)$ is a projective left $A$-module.
\end{corollary}

{\it Proof}. Let $\D_1, \ldots , \D_l$ be all the nonzero minors
from Proposition \ref{arDerA}.(1) and, for each $i=1, \ldots , l$,
$\D_i\Der_K(A) $ $\subseteq \CP_i:= \sum_k A\der_{i_1, \ldots ,
i_k;j_1, \ldots , j_r,k}\simeq A^{n-r}$ be the corresponding
inclusions. Denote by $\D$ the diagonal map $\Der_K(A)\ra \CP
:=\bigoplus_{i=1}^l\CP_i$, $\d \mapsto (\D_1\d , \ldots , \D_l\d
)$. The minors $\D_1,\ldots , \D_l$ are generators for the
Jacobian ideal $\ga_r$ which is equal to $A$ since the algebra $A$
is regular. So, $1=\sum_{i=1}^la_i\D_i$ for some elements $a_i\in
A$. The module $\CP$ is free (isomorphic to $A^{l(n-r)}$) and the
$A$-homomorphism $\g :\CP\ra \Der_K(A)$, $(u_1, \ldots ,
u_l)\mapsto \sum_{i=1}^la_iu_i$ satisfies $\g \D ={\rm
id}_{\Der_K(A)}$, the identity map. Hence the $A$-module
$\Der_K(A)$ is a direct summand of the free $A$-module $\CP$,
therefore it is a projective $A$-module. $\Box $

\begin{lemma}\label{sid}%\marginpar{sid}
 ${\rm der}_K(A)(A):=\sum_{\i \in \II_r, \j\in \JJ_{r+1}}A\derij
 (A)=\ga_r$.
\end{lemma}

{\it Proof}. Denote by $\Sigma $ the sum. By Proposition
\ref{arDerA}.(2),
 $\Sigma\subseteq \ga_r$. The nonzero minors $\D (i_1, \ldots ,
i_r;j_1, \ldots , j_r)$ form a generating set for the  ideal
$\ga_r$. Since $\der_{i_1, \ldots , i_r;j_1, \ldots , j_r,
k}(\bx_k)=\D (i_1, \ldots , i_r;j_1, \ldots , j_r)$ for any $k\in
\{1, \ldots , n\} \backslash \{ j_1, \ldots , j_r\}$, it follows
that $\ga_r\subseteq \Sigma$, hence $\Sigma =\ga_r$.
 $\Box $

A combination of Lemma \ref{sid} with Theorem \ref{7Mar05} gives a
short proof that the Jacobian ideal $\ga_r$ is invariant  for the
algebra $A$.

\begin{corollary}\label{0sid}%\marginpar{0sid}
Given another presentation, say $P_{n'}/I'$, for the algebra
$A=P_n/I$ and let $\ga_{r'}$ be the corresponding Jacobian ideal.
Then $\ga_{r'}=\s (\ga_r)$ where $\s : A\ra P_{n'}/I'$ is an
algebra isomorphism.
\end{corollary}

{\it Proof}. By Theorem \ref{7Mar05}, ${\rm der}_K(P_{n'}/I')=\s
{\rm der}_K(A)\s^{-1}$. By Lemma \ref{sid}, $\s (\ga_r)=\s ({\rm
der}_K(A)(A))=\s {\rm der}_K(A)\s^{-1}\s (A)={\rm
der}_K(P_{n'}/I')(P_{n'}/I')=\ga_{r'}$. $\Box$

The Jacobian ideal $\ga_r$ defines the {\em singular locus} of the
algebra $A$ in the sense that a prime ideal $\gp$ of $A$ contains
$\ga_r$ iff $A_\gp$ is not a regular local ring (\cite{Eisenbook},
16.20).

The next result shows that the natural derivations detect
singularities.

\begin{corollary}\label{1sid}%\marginpar{1sid}
Suppose, in addition,  that the field $K$ is  algebraically
closed, let $\gm $ be a maximal ideal of $A$. Then $\ga_r\subseteq
\gm$ $\Leftrightarrow$  ${\rm der} (A)(\gm ) \subseteq \gm$.
\end{corollary}

{\it Proof}. If $\ga_r \subseteq \gm$ then, by Lemma \ref{sid},
${\rm der} (A)(\gm ) \subseteq {\rm der} (A)(A) \subseteq
\ga_r\subseteq \gm$.

If $\ga_r \not\subseteq \gm$ then $\D (\i , \j )\not\in \gm $ for
some $\i $ and $\j =(j_1, \ldots , j_r)$ and $\der_{\i , \j ,
j_{r+1}}(\overline{x}_{j_{r+1}}-\l )=\der_{\i , \j ,
j_{r+1}}(\overline{x}_{j_{r+1}})=\D(\i , \j )\not\in \gm$ where
$\l\in K$ satisfies $\overline{x}_{j_{r+1}}-\l\in \gm$. Therefore,
${\rm der} (A)(\gm ) \not\subseteq \gm$. $\Box $

From this moment and till the end of Lemma \ref{nscoo}, we will
not use Lemma \ref{IIr1}.

For the algebra $A=P_n/I$, a set $x_{j_1}, \ldots , x_{j_r}$ is
called a {\bf non-singular} set if there exit elements, say $g_1,
\ldots , g_r\in I$, such that $\det (\overline{\frac{\der
g_i}{\der x_{j_k}}})\neq 0 $ in $A$ where
 $i,k=1, \ldots , r$. In this case we say
 that the indices $j_1, \ldots , j_r$ and the set $g_1, \ldots ,
 g_r$ are {\em non-singular}.

\begin{lemma}\label{nscoo}%\marginpar{nscoo}
Suppose that $x_{j_1}, \ldots , x_{j_r}$ is a  non-singular set,
i.e. $\det (\overline{\frac{\der g_i}{\der x_{j_k}}})\neq 0$  for
some polynomials $g_1, \ldots , g_r\in I$. Then $x_{k_1}. \ldots ,
x_{k_r}$ is a  non-singular set iff $\det (\overline{\frac{\der
g_i}{\der x_{k_\nu}}})\neq 0$.
\end{lemma}

{\it Proof}. $x_{k_1}, \ldots , x_{k_r}$ is a non-singular set iff
$\det (\overline{\frac{\der h_i}{\der x_{k_\nu}}})\neq 0$ in $A$
for some polynomials $h_1, \ldots , h_r\in I$ where $\nu =1,
\ldots , r$. The $r\times r$ matrices  $G:=(\overline{\frac{\der
g_i}{\der x_{j_\nu}}})$, $H:=(\overline{\frac{\der h_i}{\der
x_{k_\nu}}})\in M_r(Q)$ are invertible, hence the two sets of
vectors $G_i=(\overline{\frac{\der g_i}{\der x_1}},\ldots ,
\overline{\frac{\der g_i}{\der x_n}})$, $i=1, \ldots , r$ and
$H_i=(\overline{\frac{\der h_i}{\der x_1}},\ldots ,
\overline{\frac{\der h_i}{\der x_n}})$, $i=1, \ldots , r$, of the
vector space $Q^n$ (over the field $Q$) are $Q$-linear
independent, and are bases for the $r$-dimensional $Q$-vector
space $V_I\subseteq Q^n$ generated by the elements
$(\overline{\frac{\der f}{\der x_1}},\ldots , \overline{\frac{\der
f}{\der x_n}})$, $f\in I$. So, there exists an invertible matrix,
say $B\in \GL_r(Q)$, such that $B\begin{pmatrix}
G_1 \\
  \vdots \\
 G_r\\
\end{pmatrix} =\begin{pmatrix}
H_1 \\
  \vdots \\
 H_r\\
\end{pmatrix}$. Let $H':= (\overline{\frac{\der h_i}{\der
x_{j_\nu}}})$, $G':= (\overline{\frac{\der g_i}{\der
x_{k_\nu}}})\in M_r(Q)$. Then $BG'=H$ and $BG=H'$, and so the
matrices $H'$ and $G'$ are invertible, and the result follows.
$\Box $

{\bf Proof of Lemma \ref{IIr1}}. Given $\i \in \II_r$ and $\j \in
\JJ_r$, then $\D (\i , \j')\neq 0$ and $\D (\i', \j )\neq 0$ for
some $\i'$ and $\j'$. Then, by Lemma \ref{nscoo}, $\D (\i ,\j
)\neq 0$.

If $\D (\i , \j )\neq 0$ then obviously $\i \in \II_r$ and $\j\in
\JJ_r$. So, $\i \in \II_r$ and $\j \in \JJ_r$ iff $\D (\i \, \j
)\neq 0$. $\Box$

Given $\i =(i_1, \ldots , i_r)\in \II_r$ and $\j =(j_1, \ldots ,
j_r)\in \JJ_r$, then $\D =\D (\i , \j )\neq 0$. For each $j=1,
\ldots , n $, consider the vector-column
$$v_{\i ,j}\equiv v_{i_1,
\ldots , i_r; j}:= (\overline{\frac{\der f_{i_1}}{\der x_j}},
\ldots , \overline{\frac{\der f_{i_r}}{\der x_j}})^T\in
Q^r=\bigoplus_{i=1}^r Qe_i=\bigoplus_{\nu =1}^r Qv_{\i , j_\nu
}=\sum_{k=1}^n Qv_{\i , k}.$$
  Clearly, $v_{\i , j_1}\wedge
\cdots \wedge v_{\i, j_r}=\D (\i , \j) e_1\wedge \cdots \wedge
e_r$. For each $\nu =1, \ldots , r$, omitting $v_{\i , j_\nu }$ in
the exterior product, i.e. $v_{\i , j_1}\wedge \cdots \wedge
(\cdot ) \wedge \cdots \wedge v_{\i, j_r}=\Dijnu (\cdot )
e_1\wedge \cdots \wedge e_r$, gives the $Q$-linear map
$$ \Dijnu :Q^r\ra Q, \;\; q=\sum_{i=1}^rq_ie_i\mapsto \det
\begin{pmatrix}
\overline{\frac{\der f_{i_1}}{\der x_{j_1}}}& \cdots &
\overline{\frac{\der f_{i_1}}{\der x_{j_{\nu -1}}}} & q_1&
\overline{\frac{\der f_{i_1}}{\der x_{j_{\nu+1}}}}  & \cdots &
\overline{\frac{\der f_{i_1}}{\der x_{j_r}}}\\
& & & \cdots & & & \\
\overline{\frac{\der f_{i_r}}{\der x_{j_1}}}& \cdots &
\overline{\frac{\der f_{i_r}}{\der x_{j_{\nu -1}}}}& q_r&
\overline{\frac{\der f_{i_r}}{\der x_{j_{\nu +1}}}}& \cdots &
\overline{\frac{\der f_{i_r}}{\der
x_{j_r}}}\\
\end{pmatrix}
$$
Clearly, $\D^{-1}\Dijnu (\vijmu )=\d_{\nu , \mu } $ for all $\nu ,
\mu =1, \ldots , r$. The set of linear maps $\D^{-1}\D_{\i , j_\nu
}:Q^r\ra Q$, $ \nu =1, \ldots , r$ is the dual basis to the basis
$v_{\i , j_1}, \ldots , v_{\i , j_r}$ for the vector space $Q^r$.
 So, the $Q$-linear map
 $\D^{-1}\Dijnu (\cdot )\vijnu :Q^r=\bigoplus_{\mu =1}^r Q\vijmu \ra
 Q\vijnu $ is the {\em projection} onto the $\nu$'th direct
 summand, hence the identity map ${\rm id}={\rm id}_{Q^r}$ has
 the following presentation
 %\marginpar{id=Dvbi}
\begin{equation}\label{id=Dvbi}
{\rm id}=\D^{-1}\sum_{\nu =1}^r\Dijnu (\cdot ) \vijnu .
\end{equation}
If $\j'=(j_1', \ldots , j_r')$ is another non-singular set of
indices then, by Lemma \ref{nscoo}, $\D (\i , \j')\neq 0$.
$v:=(v_{\i , j_1}, \ldots , v_{\i, j_r})$ and $v':=(v_{\i , j_1'},
\ldots , v_{\i, j_r'})$ are bases for the vector space $Q^r$.
Using (\ref{id=Dvbi}), the change-of-basis matrix $(v'=vH_\i (\j ,
\j' ))$ is equal to $H_\i (\j , \j' )=\D (\i , \j )^{-1} (\Dijnu
(\vijmus ))$, the $(\nu , \mu )$-entry of it is
$$\D (\i , \j
)^{-1}\Dijnu (\vijmus ))=\D (\i , \j )^{-1}\D (\i ;j_1, \ldots ,
j_{\nu-1}, j_\mu', j_{\nu+1}, \ldots , j_r).$$
 Note that each element of the matrix $H_\i (\j ,
\j' ) $ is the ratio of two elements from the canonical set of
generators for the ideal $\ga_r$.
 Since $r$ is the
rank of the Jacobi matrix $\bJ $, the matrix $H_{\i }(\j , \j')$
{\em does not depend on} $\i $, i.e. $H_{\i }(\j , \j')=H_{\i'
}(\j , \j')$ for all non-singular $\i, \i', \j , \j'$.

Similarly, $v=v'\Hjsj  $ where $\Hjsj =\D (\i , \j')^{-1} (\Dijnus
(\vijmu ))$, the $(\nu , \mu )$-entry of the matrix is $\D (\i ,
\j')^{-1}\Dijnus (\vijmu ))=\D(\i , \j')^{-1} \D (\i , j_1',
\ldots , j_{\nu-1}', j_\mu , j_{\nu+1}', \ldots , j_r')$.

 Clearly,
$\Hjjs^{-1}=\Hjsj $ or, equivalently,
$$(\D_{\i , j_\nu'} (\vijmu
))(\Dijnu (v_{\i ,j_\mu'} ))=\D (\i , \j )\D (\i , \j')E$$ where
$E$ is the identity matrix. Equating the determinants of both
sides one gets $$\det (\D_{\i , j_\nu'} (\vijmu ))\det (\Dijnu
(v_{\i ,j_\mu'} ))=(\D (\i , \j )\D (\i , \j'))^r.$$ Clearly, $
\Hjjs H(\j', \j'')=H(\j , \j'')$ for all non-singular  $\j, \j',
\j''$.

If $\Disj \neq 0$ for some $\i'=(i_1'.\ldots , i_r')$ then $\D
(\i', \j')\neq 0$ (Lemma \ref{nscoo}). For each $i=1,\ldots , m$,
let $F_i:= (\overline{\frac{\der f_i}{\der x_1}}, \ldots ,
\overline{\frac{\der f_i}{\der x_n}})\in Q^n$. In the proof of
Lemma \ref{nscoo}, we have seen that the sets $F_{i_1}, \ldots
F_{i_r}$ and $F_{i_1'}, \ldots F_{i_r'}$ are bases for the
$Q$-vector space $V_I:=\sum_{i=1}^m QF_i\subseteq Q^n$. Let $V
(\i', \i )\in {\rm GL}_r(Q)$ be the change-of-basis matrix written
in the `vertical' form: $$\begin{pmatrix} F_{i_1'}\\ \vdots
\\F_{i_r'}\end{pmatrix} =V(\i', \i )\begin{pmatrix} F_{i_1}\\ \vdots
\\F_{i_r}\end{pmatrix} . $$

 Clearly, $V(\i', \i )^{-1}=V(\i , \i')$ and $V(\i'', \i' )V(\i' , \i
 )=V(\i'' , \i )$ for all non-singular $\i , \i', \i''$.  Let
 $$ \bJ (\i , \j )= \bJ (i_1, \ldots , i_r;j_1, \ldots , j_r):=
 (\overline{\frac{\der f_{i_k}}{\der x_{j_l}}})\in M_r(A), \;\; k,l=1, \ldots , r. $$
Then %\marginpar{J=VJH}
\begin{equation}\label{J=VJH}
\bJisjs =V(\i', \i ) \bJ (\i , \j ) \Hjjs
\end{equation}
for all non-singular $\i, \i', \j , \j'$. Since $H(\j , \j )=E$
and $V(\i , \i )=E$, it follows from (\ref{J=VJH}) that
%\marginpar{0J=VJH}
\begin{equation}\label{0J=VJH}
 V(\i', \i )=\bJ (\i', \j )\bJ (\i , \j )^{-1}\;\; {\rm and}\;\;
H(\j, \j' )=\bJ (\i , \j  )^{-1}\bJ (\i , \j' ), \end{equation}
for all non-singular $\i$, $\i'$, $\j$, \and $\j'$. Taking the
determinant of both sides of (\ref{J=VJH}), we have the equality
%\marginpar{1J=VJH}
\begin{equation}\label{1J=VJH}
\D (\i', \j') =\det (V(\i', \i )) \D (\i , \j ) \det (\Hjjs ).
\end{equation}
Finding the determinant of the matrices $V(\i' , \i )$ and $H(\j ,
\j')$ using (\ref{0J=VJH}) and then substituting them in
(\ref{1J=VJH}) we find that %\marginpar{3J=VJH}
\begin{equation}\label{3J=VJH}
\D (\i', \j')\D (\i , \j )=\D (\i', \j )\D(\i , \j').
\end{equation}

\begin{corollary}\label{2sid}%\marginpar{2sid}
For any $\i \in \II_r$ and $\j\in \JJ_r$,   $\ga_r=\sum_{\i'\in
\II_r, \j'\in \JJ_r} A\,  \det (V(\i', \i ))\, \D (\i , \j )\,
\det (H(\j, \j'))$.
\end{corollary}

{\it Proof}.  $\ga_r=\sum_{\i'\in \II_r, \j'\in \JJ_r}A\,  \D
(\i', \j')$, and the result follows from (\ref{1J=VJH}). $\Box $

Theorem \ref{25July04} describes the derivations of an {\em
arbitrary} (not necessarily regular) affine domain  (which is not
a field).

\begin{theorem}\label{25July04}%\marginpar{25July04}
 Given $\i
=(i_1, \ldots , i_r)\in \II_r$ and  $\j =(j_1, \ldots , j_r)\in
\JJ_r$, let $\{ j_{r+1}, \ldots , j_n\} =\{ 1,\ldots , n
\}\backslash \{ j_1, \ldots , j_r\}$. Then
\begin{enumerate}
\item  $\Der_K(A)=\{ \D (\i , \j )^{-1} \sum_{k=r+1}^n
a_{j_k}\der_{i_1, \ldots , i_r; j_1, \ldots , j_r,j_k}\, | $ where
the elements $ a_{j_{r+1}}, \ldots , a_{j_n}\in A$ satisfy the
following system of inclusions:
$$ \sum_{k=r+1}^n\D (\i ; j_1, \ldots , j_{\nu -1},
j_k, j_{\nu +1}, \ldots , j_r)a_{j_k}\in A\D (\i , \j ), \;\; \nu
=1, \ldots , r\}. $$
 \item $\Der_K(A)=\{ \sum_{k=r+1}^na_{j_k}\der_{j_k}-
 \D (\i , \j )^{-1}\sum_{\nu =1}^r \Dijnu
 (\sum_{k=r+1}^na_{j_k}\vijk )\der_{j_\nu}\,| $ where
the elements $ a_{j_{r+1}}, \ldots , a_{j_n}\in A$ satisfy the
following system of inclusions: $$\D (\i , \j )^{-1}\sum_{\nu
=1}^r \Dijnu
 (\sum_{k=r+1}^na_{j_k}\vijk )\in A, \;\; \nu=1,
 \ldots , r\}.$$ The geometric meaning of these inclusions (in both statements) is that
 all the coordinates of the vector $ \sum_{k=r+1}^na_{j_k}\vijk
 \in \bigoplus_{\nu =1}^r Q\vijnu $ in the basis $v_{\i , j_1},
 \ldots , v_{\i , j_r}$ belong to $A$.
\end{enumerate}
\end{theorem}

{\it Proof}.  Let $\D =\D (\i , \j )$ and  $\d \in \Der_K(A)$.
Then
 $\D \d =\sum_{k=r+1}^n
a_{j_k}\der_{i_1, \ldots , i_r; j_1, \ldots , j_r,j_k}$ (by
Proposition \ref{arDerA}) for  unique elements $a_{j_k}\in A$. By
(\ref{1dbibj}),
\begin{eqnarray*}
\d & = & \D^{-1} (\D \sum_{k=r+1}^na_{j_k}\der_{j_k}-\sum_{\nu
=1}^r (\sum_{k=r+1}^n\D (\i ; j_1, \ldots , j_{\nu -1}, j_k,
j_{\nu+1}, \ldots , j_r)a_{j_k})\der_{j_\nu })\\
&=& \sum_{k=r+1}^na_{j_k}\der_{j_k}-\sum_{\nu =1}^r
(\D^{-1}\sum_{k=r+1}^n\D (\i ; j_1, \ldots , j_{\nu -1}, j_k,
j_{\nu+1}, \ldots , j_r)a_{j_k})\der_{j_\nu }\\
&=&\sum_{k=r+1}^na_{j_k}\der_{j_k}-\D (\i , \j )^{-1}\sum_{\nu
=1}^r \Dijnu
 (\sum_{k=r+1}^na_{j_k}\vijk )\der_{j_\nu}.
\end{eqnarray*}
So, both statements are equivalent. Since $\d \in \Der_K(A)$, $\d
(\bx_{j_\nu })\in A$ for all $\nu =1, \ldots , r$. This proves the
necessary conditions.

Conversely, suppose that the conditions of the theorem are
satisfied for a derivation $\d =\D (\i , \j )^{-1} \sum_{k=r+1}^n
a_{j_k}\der_{i_1, \ldots , i_r; j_1, \ldots , j_r,j_k}$. Then
$\d\in \Der_K(A)$ since $\d (\bx_i)\in A$ for $i=1,\ldots , n$.
$\Box $

{\it Remark}. By (\ref{0J=VJH}), the system of inclusions in
Theorem \ref{25July04} is equivalent to the system of inclusions
%\marginpar{dHjn}
\begin{equation}\label{dHjn}
\sum_{k=r+1}^n \det (H(\j ; j_1, \ldots , j_{\nu-1}, j_k,
j_{\nu+1}, \ldots , j_r))a_{j_k}\in A, \; \nu =1, \ldots , r,
\end{equation} and it does not depend on $\i $. So, any $r$
algebraically independent elements of the ideal $I$ give the same
coefficients of the above system. So,  any $r$ algebraically
independent elements of the ideal $I$ `determine' (in this way)
the module of derivations $\Der_K(A)$.

\begin{corollary}\label{c425July04}%\marginpar{c425July04}
$\Der_K(A)_\gp \simeq A^{n-r}_\gp$ for each prime ideal $\gp $ of
$A$ such that $\ga_r\not\subseteq \gp$.
\end{corollary}

{\it Proof}. Since $\ga_r\not\subseteq \gp$, there exists $\D =
\D( \i , \j )\not\in \gp$, and so $\D $ is a  unit of the local
ring $A_\gp$. Since $\Der_K(A)_\gp \simeq \Der_K(A_\gp )$, then by
Theorem \ref{25July04}.(i),
$$ \Der_K(A_\gp )=\bigoplus_{k=r+1}^nA_\gp \der_{i_1, \ldots , i_r;
j_1, \ldots , j_r,j_k}\simeq A_\gp^{n-r}. \;\; \Box $$

\begin{corollary}\label{c025July04}%\marginpar{c025July04}
The rank of the $r\times (n-r)$ matrix $(\D (\i , j_1, \ldots ,
j_{\nu -1}, j_k, j_{\nu +1}, \ldots , j_r))$, $\nu =1, \ldots ,
r$; $k=r+1, \ldots , n$,  in Theorem \ref{25July04}.(1) is equal
to $\dim_Q(\sum_{k=r+1}^nQ\vijk )=$ the rank of the $r\times
(n-r)$ matrix $(\overline{\frac{\der f_i}{\der x_{j_k}}})$,
$i=i_1, \ldots , i_r$; $k=r+1, \ldots , n$.
\end{corollary}

{\it Proof}. By (\ref{id=Dvbi}), the first matrix  multiplied by
$\D^{-1}$ (resp. the second matrix) in the corollary is the matrix
of coordinates of the vectors $\vijk$,$k=r+1, \ldots , n$ in the
basis $v_{\i , j_1}, \ldots , v_{\i , j_r}$ (resp. $e_1, \ldots ,
e_r$). Now, the result is evident. $\Box $

Clearly, a derivation $\d $ of $A$ is uniquely determined by its
action on generators of the algebra $A$. For the polynomial
algebra $P_n$ one has a nice formula $\d =\sum_{i=1}^n\d
(x_i)\der_i$. The next corollary gives a similar {\em explicit
formula} for $\d $ in general situation.

\begin{corollary}\label{c25July04}%\marginpar{c25July04}
Any derivation $\d \in \Der_K(A)$ of the algebra $A$  is uniquely
determined by its action on any subset of the canonical
generators,  say $\bx_{j_{r+1}}, \ldots , \bx_{j_n}$, that are
complement to a non-singular subset of the canonical generators,
say $\bx_{j_1}, \ldots , \bx_{j_r}$:
$$ \d =\sum_{k=r+1}^n \d (\bx_{j_k})\der_{j_k}-\D(\i , \j
)^{-1}\sum_{\nu =1}^r\Dijnu (\sum_{k=r+1}^n\d (\bx_{j_k})\vijk
)\der_{j_\nu }$$ where $\i =(i_1, \ldots , i_r)$ is any element of
$\II_r$.
\end{corollary}

{\it Proof}. By Theorem \ref{25July04}, $\d
=\sum_{k=r+1}^na_{j_k}\der_{j_k}-
 \D (\i , \j )^{-1}\sum_{\nu =1}^r \Dijnu
 (\sum_{k=r+1}^na_{j_k}\vijk )\der_{j_\nu}$ for unique elements
 $a_{j_k}\in A$. Clearly, $a_{j_k}=\d (\bx_{j_k})$, and the result
 follows. $\Box $

\begin{corollary}\label{c225July04}%\marginpar{c225July04}
Keep the notation of Theorem \ref{25July04}. Given   $\i , \i'\in
\II_r$ and $ \j\in\JJ_r$, then for  each $\nu =1, \ldots , r$:
$$ \D (\i , \j )^{-1} \Dijnu (\sum_{k=r+1}^na_{j_k}\vijk )=
\D (\i' , \j )^{-1} \D_{\i', j_\nu } (\sum_{k=r+1}^na_{j_k}v_{\i',
j_k})$$
 for all $a_{j_{r+1}}, \ldots , a_{j_n}\in Q$.
\end{corollary}

{\it Proof}. By linearity, it suffices to prove the equality for
all elements $a_{j_{r+1}}, \ldots , a_{j_n}\in \D (\i , \j )\D
(\i', \j )A$. For any such a choice of the elements, they
determine uniquely a derivation $\d $ from Corollary
\ref{c25July04} where $\d (\bx_{j_\nu} )=a_{j_\nu}$, $\nu = 1,
\ldots , r$ (since they satisfy the conditions of Theorem
\ref{25July04}.(2)). We have the same derivation $\d $ for each
pair $(\i ,\j)$ and $(\i', \j )$. Now, comparing the coefficients
of $\der_{j_\nu }$ in two expressions for  $\d $ in Corollary
\ref{c25July04} we finish the proof (see (\ref{DerPnI})).  $\Box $

For a maximal ideal $\gm $ of $A$, $k_\gm:=A_\gm /\gm A_\gm \simeq
A/\gm$ is the {\em residue} field for $\gm$. Note that $[k_\gm
:K]<\infty$, so $k_\gm$ is a {\em finite separable} extension of
the field $K$. $\Der_K(A_\gm , k_\gm )$ is called the {\em tangent
space} at the point $x=x_\gm $ of the algebraic variety $V=\Spec
(A)$ denoted $T(V)_x$. There is the  canonical  isomorphism of
left $k_\gm $-modules %\marginpar{TVx}
\begin{equation}\label{TVx}
\Der_K(A_\gm , k_\gm )\ra (\gm / \gm^2)^*:= \Hom_{k_\gm }(\gm /
\gm^2, k_\gm ), \;\; \d \mapsto (\overline{\d}: a+\gm^2\mapsto \d
(a) +\gm ).
\end{equation}
This fact is obvious  if the field $K$ is an algebraically closed
field since then $A/\gm^2= K\oplus \gm/\gm^2$. If the field $K$ is
not necessarily  algebraically closed  then using the {\em
Hensel's Lemma} and the fact that the ring $A/\gm^2$ is complete
in the $\gm$-adic topology and $k_\gm $ is separable over $K$ one
can find an isomorphic copy of the field $k_\gm $ in $A/\gm^2$ so
that $A/\gm^2\simeq k_\gm \oplus \gm / \gm^2$. Since $\d (k_\gm
)=0$ as $k_\gm $ is separable over $K$, the result follows. The
canonical epimorphism $A_\gm \ra k_\gm $ yields the canonical
$k_\gm $-homomorphism $\Der_K(A_\gm )\ra \Der_K(A_\gm , k_\gm )$.
In combination with (\ref{TVx}), we have the canonical $k_\gm
$-module homomorphism

%\marginpar{1TVx}
\begin{equation}\label{1TVx}
\frac{\Der_K(A_\gm )}{\gm \Der_K(A_\gm )}\ra (\gm / \gm^2)^*, \;\;
\d +\gm \Der_K(A_\gm  )\mapsto (\overline{\d}:a+\gm^2\mapsto \d
(a)+\gm ).
\end{equation}

\begin{corollary}\label{c325July04}%\marginpar{c325July04}
For any maximal ideal $\gm$ of $A$ such that $\ga_r\not\subseteq
\gm$, the map (\ref{1TVx}) is an isomorphism (so, at a simple
point $x=x_\gm$ of the variety $V$ the tangent space $T(V)_x$ is
canonically isomorphic to $\Der_K(A_\gm )/\gm \Der_K(A_\gm )$).
\end{corollary}

{\it Proof}. Since $\ga_r\not\subseteq \gm$, there exists $\D =
\D( \i , \j )\not\in \gm$, and so $\D $ is a  unit of the local
ring $A_\gm$. Since $\Der_K(A)_\gm \simeq \Der_K(A_\gm )$, then,
by Theorem \ref{25July04}.(1),
$$ \Der_K(A_\gm )=\bigoplus_{k=r+1}^nA_\gm \der_{i_1, \ldots , i_r;
j_1, \ldots , j_r,j_k}\simeq A_\gm^{n-r},  $$ and so $T':=
\Der_K(A_\gm )/\gm \Der_K(A_\gm )\simeq k_\gm^{n-r}$, hence
$\dim_{k_\gm} (T')=n-r=\Kdim (A_\gm )=\dim_{k_\gm }(\gm/\gm^2)^*$
since the point $x=x_\gm$ is simple. Since $\D^{-1}\der_{\i ; \j,
j_k}(\bx_{j_l}+\gm^2)=\d_{k,l}+\gm$ (the Kronecker delta) where
$k,l=r+1, \ldots , n$, the map $T'\ra (\gm / \gm^2)^*$ is
injective, and so $T'\simeq (\gm / \gm^2)^*$ as the vector spaces
have the same dimension over $k_\gm$. $\Box$

{\it Definition}. For each maximal ideal $\gm $ of $A$ such that
$\ga_r\subseteq \gm$, let $\CK_\gm $, $\CI_\gm $, and $\CC_\gm $
be the {\em kernel}, the {\em image}, and the {\em cokernel} of
the $k_\gm$-homomorphism (\ref{1TVx}), and we have the short exact
sequence of $k_\gm$-modules: %\marginpar{CKCm}
\begin{equation}\label{CKCm}
0\ra \CK_\gm \ra \frac{\Der_K(A_\gm )}{\gm \Der_K(A_\gm )}\ra (\gm
/ \gm^2)^*\ra \CC_\gm \ra 0.
\end{equation}
 The vector spaces involved in
the short exact sequence and
 $\CI_\gm $ are important invariants of the singularity at
 $x=x_\gm$.

For each $s$ such that $1\leq s\leq r$, let
$$ \gder_K(A)_s:=\{ \d \in \Der_K(A)\, | \, \d (A)\subseteq
\ga_s\}.$$ Then
$${\rm der}_K(A)\subseteq \gder_K(A)_r\subseteq\gder_K(A)_{r-1}\subseteq
\cdots \subseteq \gder_K(A)_1\subseteq \Der_K(A)$$ is the
ascending chain of left $A$-modules such that
$$[\gder_K(A)_s,\gder_K(A)_s]\subseteq \gder_K(A)_s$$ for each
$s$. So, each $\gder_K(A)_s$ is a Lie subalgebra of the Lie
algebra $\Der_K(A)$. It follows immediately that each Lie algebra
$ \gder_K(A)_s$ is an {\em invariant} for the algebra $A$ (i.e.
$\s \gder_K(P_n/I)_s\s^{-1}=\gder_K(P_{n'}/I')_s$).

For each $s$ such that $1\leq s\leq r$, let $\i=(i_1, \ldots ,
i_s)$ and $\j =(j_1, \ldots , j_{s+1})$ satisfy $1\leq i_1<\cdots
<i_s\leq m$ and $1\leq j_1<\cdots <j_{s+1}\leq n$. Then
%\marginpar{sdbibj}
\begin{eqnarray}\label{sdbibj}
 \derij  =  \der_{i_1, \ldots , i_s; j_1, \ldots , j_{s+1}}:= \det
 \begin{pmatrix}
  \overline{\frac{\der f_{i_1}}{\der x_{j_1}}} & \cdots &  \overline{\frac{\der f_{i_1}}{\der
  x_{j_{s+1}}}}\\
  \vdots & \vdots & \vdots \\
  \overline{\frac{\der f_{i_s}}{\der x_{j_1}}} & \cdots &  \overline{\frac{\der f_{i_s}}{\der
  x_{j_{s+1}}}}\\
  \der_{j_1}& \cdots & \der_{j_{s+1}}\\
\end{pmatrix}\in \Der_K(A, A/\ga_{s+1})
\end{eqnarray}
where the bar means mod $\ga_{s+1}$ and $\ga_{r+1}:=0$.

{\it Definition}. For the algebra $A=P_n/I$ and a given set $f_1,
\ldots , f_m$ of generators for the ideal $I$, and for each $s$
such that $1\leq s\leq r$, we denote by ${\rm der}_K(A)_s$ the
left $A/\ga_{s+1}$-submodule of $\Der_K(A, A/\ga_{s+1})$ generated
by all the derivations $\derij$ from (\ref{sdbibj}).

The derivations $\derij$ {\em depend} on the choice of generators
for the ideal $I$. The next theorem shows that the left
$A/\ga_{s+1}$-module ${\rm der}_K(A)_s$ does not depend on the
choice of generators for the ideal $I$ and on the choice of
presentation of the algebra $A$ as the factor algebra $P_n/I$. We
call ${\rm der}_K(A)_s$ the  $A/\ga_{s+1}$-module of {\bf natural
derivations} from the algebra $A$ to $A/\ga_{s+1}$. So, higher
natural derivations, ${\rm der}_K(A)_s$, $1\leq s\leq r$, are {\em
invariants} for the algebra $A$.

\begin{theorem}\label{8Mar05}%\marginpar{8Mar05}
For each $s$ such that $1\leq s\leq r$,
\begin{enumerate}
\item the $A/\ga_{s+1}$-module ${\rm der}_K(A)_s$ does not depend
on the choice of generators for the ideal $I$, and  \item  the
$A/\ga_{s+1}$-module ${\rm der}_K(A)_s$ does not depend on the
presentation of the algebra $A$ as the factor algebra $P_n/I$.
\end{enumerate}
\end{theorem}

{\it Proof}.  $1$.  Given another set, say $g_1, \ldots , g_l$, of
generators for the ideal $I$. Let ${\rm der}_K(A)_s'$ be the
corresponding $A/\ga_{s+1}$-submodule of $\Der_K(A,A/\ga_{s+1})$.
Then $f_i=\sum_{k=1}^la_{ik}g_k$ for some $a_{ik}\in P_n$.
Repeating the arguments of the proof of Theorem \ref{7Mar05}.(1),
we have, for the derivation $\derij$ from (\ref{sdbibj}),  the
equality %\marginpar{sdijdijp}
\begin{equation}\label{sdijdijp}
\der_{\i , \j}=\sum_{\i'}\det (a_{\i, \i'})\der_{\i', \j}'
\end{equation}
where $\i'=(i_1', \ldots , i_s')$ runs through all the $s$-tuples
such that $1\leq i_1'<\cdots <i_s'\leq l$, $a_{\i, \i'}$ is an
$s\times s$ matrix with $(\nu, \mu )$-entry equal to
$\overline{a_{i_\nu , i_{\mu}'}}\in A/\ga_{s+1}$, and $\der_{\i',
\j}'$ is the derivation (\ref{sdbibj}) for the second choice of
generators. It follows that ${\rm der}_K(A)_s\subseteq {\rm
der}_K(A)_s'$. By symmetry, the opposite inclusion is also true
which proves that ${\rm der}_K(A)_s={\rm der}_K(A)_s'$.

$2$. The proof of this statement is similar to the proof of the
second statement of Theorem \ref{7Mar05}, wee keep the notation
from there. Given another presentation $P_{n'}/I'$ for the algebra
$A$ where $P_{n'}=K[x_1', \ldots , x_{n'}']$, etc. Then the
canonical generators for ${\rm der}_K(P_{n+n'}/I'')_s$ (see the
proof of Theorem \ref{7Mar05}.(2)) has the form
$\derij'':=\der_{\i, n+1, \ldots, n+n'; \j , n+1, \ldots, n+n'}$
where $\i =(i_1, \ldots , i_s)$ and $\j =(j_1, \ldots , j_{s+1})$.
The actions of the derivations $\derij''$ and $\derij$ on the
generators $\bx_1, \ldots , \bx_n$ of the algebra
$A=P_n/I=P_{n+n'}/I''$ coincide, therefore $\derij''=\derij$, and
so ${\rm der}_K(A)_s\simeq {\rm der}_K(A)_s''$. By symmetry, we
have ${\rm der}_K(A)_s'\simeq {\rm der}_K(A)_s''$, and the result
follows, where ${\rm der}_K(A)_s'$ and ${\rm der}_K(A)_s''$ are
the natural derivations for the presentations $P_{n'}/I'$ and
$P_{n+n'}/I''$ respectively.  $\Box $

For each $s$ such that $1\leq s\leq r$, let
$\overline{\gder}_K(A)_s$ be the image of the natural map
$$ \gder_K(A)_s\ra \Der_K(A, A/\ga_{s+1}), \d \mapsto (a\mapsto \d
(a)+\ga_{s+1}).$$ Clearly, $\overline{\gder}_K(A)_s$ is a left
$A/\ga_{s+1}$-submodule of $\Der_K(A, A/\ga_{s+1})$ and ${\rm
der}_K(A)_s\subseteq  \overline{\gder}_K(A)_s$.

{\it Question}. {\em When}  ${\rm der}_K(A)_s =
\overline{\gder}_K(A)_s$?

\begin{corollary}\label{c8Mar05}%\marginpar{c8Mar05}
For each $s$ such that $1\leq s\leq r$,
\begin{enumerate}
\item ${\rm
der}_K(A)_s(A)=\overline{\gder}_K(A)_s(A)=\ga_s/\ga_{s+1}$, \item
${\rm der}_K(A)_s$ and $\overline{\gder}_K(A)_s$ are invariants
under the action of the group $\Aut_K(A)$.
\end{enumerate}
\end{corollary}

{\it Proof}. $1$. Since $\der_{i_1, \ldots , i_s;j_1, \ldots ,
j_s, k}(\bx_k)=\D (i_1, \ldots , i_s;j_1, \ldots , j_s)+\ga_{s+1}$
we have $\ga_s/\ga_{s+1}={\rm der}_K(A)_s(A)\subseteq
\overline{\gder}_K(A)_s(A)\subseteq \ga_s/\ga_{s+1}$.

$2$. Since  the Jacobian ideals are invariants under the action of
the group $\Aut_K(A)$, the result for ${\rm der}_K(A)_s$ follows
from Theorem \ref{8Mar05} and for $\overline{\gder}_K(A)_s$ from
the definition.  $\Box $

%%%%%%%%%%%%%%%%%% SECTION 3 %%%%%%%%%%%%%%%%%%%%%%%%

\section{Defining relations for the ring of differential operators
on a smooth irreducible affine algebraic
variety}\label{Reldefops}%\marginpar{Reldefops}

In this section, Theorems \ref{9bFeb05}, \ref{9Feb05} and
\ref{25Feb05} are proved. On the way, a short direct proof is
given of the fact that the algebra of differential operators on a
regular affine domain is a simple Noetherian algebra generated by
the affine domain and its derivations (Theorem \ref{3Feb05}).

\begin{theorem}\label{CDAffl}%\marginpar{CDAffl}
Let the algebra $A$ be a regular algebra. Then $\CD (A)\ra
\prod_{\i\in \II_r , \j\in \JJ_r}\CD (A)_{\D (\i , \j )}$ is a
left and right faithfully flat extension of algebras where $\CD
(A)_{\D (\i , \j )}$ is the localization of the algebra $\CD (A)$
at the powers of the element $\D (\i , \j )$.
\end{theorem}

{\it Proof}. The algebra $A$ is regular, so $A=\ga_r=(\D (\i , \j
))$, hence the ideal of $A$ generated by any power of the elements
$\D_{\i , \j}$ is also equal to $A$. The extension is a flat
monomorphism. Suppose that the extension is not, say left
faithful, then there exists a proper left ideal, say $L$, of
$\CD(A)$ such that $\prod_{\i , \j}\CD (A)_{\D (\i , \j )}\t_{\CD
(A)}(\CD (A)/L)=0$, equivalently, there exists a sufficiently
large natural number $k$ such that $\D (\i , \j )^k\in L$ for all
$\i , \j $. Since $A=(\D (\i , \j )^k)\subseteq L$, we must have
$L=\CD (A)$, a contradiction. $\Box $

Let $R$ be a (not necessarily commutative)  algebra over a field
$K$, and let $\d $ be a
 $K$-derivation of the algebra $R$. For any elements $a,b\in R$
 and a natural number $n$, an easy induction argument gives the
 equality
 $$ \d^n(ab)=\sum_{i=0}^n\, {n\choose i}\d^i(a)\d^{n-i}(b).$$
 It follows that the kernel $C(\d , R):=\ker \, \d $ of $\d $ is a
 subalgebra (of constants for $\d $) of $R$ (since $\d (ab)=\d
 (a)b+a\d (b)=0$ for $a,b\in C(\d , R)$), and the union of the
 vector spaces $N(\d ,R)=\cup_{i\geq 0}\, N(\d , i,R)$ is a positively
 filtered algebra (so-called, the {\em nil-algebra} of $\d $) where $N(\d
 , i,R):=\{ a\in R\, | \, \d^{i+1}(a)=0\}$,  that is
 $$N(\d , i,R)N(\d , j,R)\subseteq N(\d , i+j,R), \;\; {\rm for \;\;
 all}\;\; i,j\geq 0.$$
 Clearly, $N(\d , 0,R)= C(\d , R)$ and $N(\d , R):=\{ a\in R \, | \ \d^n (a)=0$
  for some natural $n\}$.

A $K$-derivation $\d $ of
 the algebra $R$ is a {\em locally nilpotent } derivation if for
 each element $a\in R$ there exists a natural number $n$ such
 that $\d^n(a)=0$. A $K$-derivation $\d $ is locally nilpotent iff
 $R=N(\d , R)$. A derivation of $R$ of the type $\ad (r) :x\mapsto
 [r,x]:=rx-xr$ is called an {\em inner} derivation of $R$ where
 $r\in R$.

The {\em Zariski Lemma} is instrumental in dealing with
derivations on complete local commutative rings.

\begin{lemma}\label{zarlemma}%\marginpar{zarlemma}
({\bf Zariski Lemma}, \cite{ZarAMJ1965}) Let $C$ be a complete
local  commutative ring and $d$ be a derivation of $C$ sending a
non-unit into a unit. Then $C=C_0[[x]]$ for some ring $C_0$ with
$d(C_0)=0$.
\end{lemma}

We need a {\em noncommutative} analogue of the Zariski Lemma
(Lemma \ref{dx=1}).

Given a ring $R$ and its derivation $d$. The {\em Ore extension}
$R[x;d]$ of $R$ is a ring freely generated over $R$ by $x$ subject
to the defining relations: $xr=rx+d(r)$ for all $r\in R$.
$R[x;d]=\bigoplus_{i\geq 0}Rx^i=\bigoplus_{i\geq 0}x^iR$ is a left
and right free $R$-module. Iterating the construction  one gets an
{\em iterated Ore extension} $R[x_1; d_1]\cdots [x_n; d_n]$. A
particular kind of an iterated Ore extension will appear in the
proof of Theorem \ref{CDAd=gd}: given a set {\em commuting}
derivations $d_1, \ldots , d_n$ of the ring $R$, and a set
commuting indeterminates $t_1, \ldots , t_n$, then $S:= R[t_1,
\ldots , t_n; d_1,\ldots , d_n]$ is a ring freely generated by $R$
and (commuting) elements $t_1, \ldots , t_n$ subject to the
defining relations: $t_ir=rt_i+d_i(r)$ for all $r\in R$ and $i=1,
\ldots , n$. $S=\bigoplus_{\alpha \in \mathbb{N}^n} Rt^\alpha =
\bigoplus_{\alpha \in \mathbb{N}^n} t^\alpha R$ is a left and
right free $R$-module. The ring $S$ is Noetherian iff the ring $R$
is so.

\begin{lemma}\label{dx=1}%\marginpar{dx=1}
 Let $R$ be an algebra over a field $K$ of characteristic zero and
$\d $ be a $K$-derivation of $R$ such that $\d (x)=1$ for some
$x\in R$. Then $N(\d ,R)=C(\d , R)[x; d]$ is the Ore extension
with coefficients from the subring $C(\d , R):=\ker \, \d $ of
constants of the derivation $\d $, and the derivation $d$ of the
algebra $ C(\d , R)$ is the restriction of the inner derivation
$\ad \, x $ of the algebra $R$ to its subalgebra $C(\d , R)$. For
each $n\geq 0$, $N(\d , n, R)=\bigoplus_{i=0}^n\, C(\d , R)x^i$.
\end{lemma}

{\it Proof}. For each element $c\in C:=C(\d , R)$,
$$ \d ([x,c])=[\d (x), c]+[x, \d (c)]=[1,c]+[x,0]=0,$$
thus $d (C)\subseteq C$, and $d $ is a $K$-derivation of the
algebra $C$.

First, we show that the $K$-subalgebra $N'$ of $N:=N(\d , R)$
generated by $C$ and $x$ is the Ore extension $C[x;d]$. We have
$N'=\sum_{i\geq 0}\, Cx^i$ since, for each $c\in C$, $xc-cx=d(c)$.
So, it remains to prove that the sum $\sum_{i\geq 0}\, Cx^i$ of
left $C$-modules is a direct sum. Suppose this is not the case,
then there is a nontrivial relation of degree $n>0$,
$$ c_0+c_1x+\cdots +c_nx^n=0, \;\; c_i\in C,\;\; c_n\neq 0.$$
We may assume that the degree $n$ of the relation above is the
least one. Then applying $\d $ to the relation above we obtain the
relation
$$ c_1+2c_2x+\cdots +nc_nx^{n-1}=0$$
of smaller degree $n-1$ since $nc_n\neq 0$ (char $K=0$), a
contradiction. So, $N'=C[x;d]$.

It remains to prove  that $N=N'$. The inclusion $N'\subseteq N$ is
obvious. In order to prove the inverse inclusion it suffices to
show that all subspaces $N(\d , i, R)$ belong to $N'$. We use
induction on $i$. The base of the induction is trivial since $N(\d
, 0, R)=C$. Suppose that $i>0$, and $N(\d , i-1, R)\subseteq N'$.
Let $u$ be an arbitrary element of $N(\d, i, R)$. Then $\d (u)\in
N(\d , i-1, R)\subseteq N'$. For an arbitrary element $a=\sum \,
c_jx^j\in N'$, we have $\d (b)=a$ where $b=\sum\,
(j+1)^{-1}c_jx^{j+1}\in N'$. Therefore, in the case of $a=\d
(u)\in N'$, we have $\d (u)=\d (b)$ for some $b\in N'$. Hence, $\d
(u-b)=0$, and $u\in b+C\subseteq N'$. This means that $N=N'$, as
required. $\Box $

\begin{theorem}\label{CDAd=gd}%\marginpar{CDAd=gd}
Given  $\i =(i_1, \ldots , i_r) \in \II_r$ and $\j =(j_1, \ldots ,
j_r)\in \JJ_r$, i.e. $\D =\D (\i , \j )\neq 0$, and  $\{ j_{r+1},
\ldots j_n\}=\{ 1, \ldots , n\} \backslash \{ j_1, \ldots , j_r\}$
and let  $A_\D$ be the localization of the algebra $A$ at the
powers of the element $\D$. Then
\begin{enumerate}
\item the algebra $\CD (A_\D )$ of differential operators on
$A_\D$ is a simple Noetherian algebra equal to the algebra  $ A_\D
\langle \der_{\i ; \j, j_{r+1}},\ldots ,\der_{\i ; \j,
j_{n}}\rangle $ (generated by $A_\D $ and $\Der_K(A_\D )$ in ${\rm
End}_K(A_\D )$). \item $\Der_K(A_\D )=\bigoplus_{\nu =r+1}^nA_\D
\der_{\i ; \j, j_{\nu}}$.
\end{enumerate}
\end{theorem}

{\it Proof}. Without loss of generality we can assume that $\i
=(1,2, \ldots , r)$ and $\j =(1,2, \ldots , r)$. Let
$$\der_{r+1}:=\D^{-1}\der_{\i ; \j, r+1},\ldots  ,\der_n:=\D^{-1}\der_{\i ; \j,
n}. $$ Then $\der_i(x_j)=\d_{ij}$ for all $i,j=r+1, \ldots n $.
Each derivation $[\der_i,\der_j]$ (the commutator of derivations)
annihilates the polynomial subalgebra $K[x_{r+1}, \ldots , x_n]$
of the finitely generated domain $A$, and both algebras have the
same Krull dimension $n-r$. Since $K$ has characteristic zero, it
follows that the derivations $\der_i$ commute. Clearly, the inner
derivations  $\ad (x_{r+1}), \ldots , \ad (x_n)$ of the algebra
$E:={\rm End}_K(A_\D )$ commute. Applying several times Lemma
\ref{dx=1} we see that the algebra
$$ N=N(\ad (x_{r+1}), \ldots , \ad
(x_n); E):= \bigcap_{i=r+1}^nN(\ad (x_i), E)=C[\der_{r+1}, \ldots
, \der_n; \ad (\der_{r+1}), \ldots , \ad (\der_n)]$$ is an
iterated Ore extension with coefficients from the algebra
$C:=\cap_{i=r+1}^n\ker (\ad (x_i))$. So, any element $u$ of $N$ is
uniquely written as a sum $u=\sum_{\alpha \in
\mathbb{N}^{n-r}}c_\alpha \der^\alpha$, $ c_\alpha \in C$. The
algebra $N=\cup_{i\geq 0}N_i$ has a natural filtration by the
total degree of the $\der_i$th. Clearly, $\CD (A_\D )\subseteq N$
and $\CD (A_\D )_i\subseteq N_i$ for each $i\geq 0$. Let us prove,
by induction on $i$, that $\CD (A_\D )_i=D_i:=\sum_{|\alpha|\leq
i}A_\D \der^\alpha$. The case $i=0$ is true, $\CD (A_\D )=A_\D
=D_0$. Suppose that  $i>0$, and by induction $\CD (A_\D
)_{i-1}=D_{i-1}$. Take $u\in \CD (A_\D )_i$. Since $\CD (A_\D
)_i\subseteq N_i$, $u=\sum_{|\alpha |\leq i}c_\alpha \der^\alpha$
for some $c_\alpha \in C$. For each $j=r+1, \ldots , n$, $\ad
(x_j)(u)=\sum_{|\alpha |\leq i}\alpha_jc_\alpha \der^{\alpha
-e_j}\in \CD(A_\D )_{i-1}=D_{i-1}$, therefore all $c_\alpha \in
A_\D$ such that $\alpha \neq 0$. Since $c_0=u-\sum_{\alpha \neq 0,
| \alpha | \leq i} c_\alpha \der^\alpha\in C\cap \CD (A_\D )_i$,
it follows from the claim below that $c_0\in A_\D $. Therefore
$\CD (A_\D )_i=D_i$, and so $\CD (A_\D )=A_\D \langle \der_{r+1},
\ldots ,\der_n\rangle$. Since, on the one hand, $\CD (A_\D
)_1=A_\D \oplus \Der_K(A_\D )$ (this is true for any commutative
algebra), and, on the other hand, $\CD (A_\D )_1=A_\D \bigoplus
(\sum_{i=r+1}^nA_\D \der_i)$, we must have $\Der_K(A_\D
)=\sum_{i=r+1}^nA_\D \der_i$. So, the algebra $\CD (A_\D )$ is
generated by the algebra $A_\D$ and the derivations $\Der_K(A_\D
)$. The derivations $\der_{r+1}:=\D^{-1}\der_{\i ; \j, r+1},\ldots
,\der_n:=\D^{-1}\der_{\i ; \j, n} $ are equal to the partial
derivatives $\frac{\der}{\der x_{r+1}}, \ldots , \frac{\der}{\der
x_n}$ respectively. So, the algebra $\CD (A_\D )$ is a factor
algebra of the iterated Ore extension $A_\D [t_{r+1}, \ldots ,
t_n; \frac{\der}{\der x_{r+1}}, \ldots , \frac{\der}{\der x_n}]$
which is a Noetherian algebra since the algebra $A_\D $ is so.
Then $\CD (A_\D )$ is a Noetherian algebra.

{\it Claim}. $C\cap \CD (A_\D )_i=A_\D $ {\em for all } $i\geq 0$.

We use induction on $i$. The case  when $i=0$ is trivial, $C\cap
\CD (A_\D )_0=C\cap A_\D =A_\D $. Note that  $C$ is an $A_\D
$-bimodule, and so is invariant under $\ad (a)$ for any $a \in
A_\D $. If the intersection $I_i:= C\cap \CD (A_\D )_i\neq A_\D$
for some $i\geq 1$, then obviously $I_1\neq A_\D $. Fix an
element, say $c\in I_1\backslash A_\D $,  and an element $a\in
A_\D $ such that $0\neq b:=\ad (a) (c)\in I_0=A_\D $. The map
$c:A_\D \ra A_\D $ is a $K[x]$-module homomorphism where
$K[x]=K[x_{r+1}, \ldots x_n]$. Let $Q$ and $K(x)$ be the field of
fractions of the algebras $A_\D $ and $K[x]$ respectively. The
fields $Q$ and $K(x)$ are finitely generated over $K$ and have the
same transcendence degree $n-r$ over $K$. Therefore
$[Q:K(x)]<\infty$ and $Q=K(x)\t_{K[x]}A_\D$. Localizing at
$K[x]\backslash \{ 0\}$, one can extend uniquely the map   $c$ to
a $K(x)$-linear map from $Q$ to $Q$.  Let $p(t)\in K(x)[t]$ be a
monic polynomial of least degree such that $p(c)=0$. Then $0=(\ad
\, a) (p(c))=p'(a)b$ where $p'=\frac{dp}{dx}$, and so $p'(a)=0$.
This contradicts the minimality of $p$ (char $K=0$). Therefore,
$I_i$ must be equal to $A_\D$ for all $i\geq 0$.

Let $L$ be a nonzero ideal of $\CD (A_\D )$. It remains to prove
that $L=\CD (A_\D )$. Take a nonzero element, say $u$, of $L$.
Applying several times maps of the type $\ad (x_{j_k})$, $r+1\leq
k\leq n$, to the element $u$ we have a nonzero element, say
$u_1\in L\cap A_\D$. Since $\Kdim (A_\D )=\Kdim (K[x_{r+1}, \ldots
, x_n])$, we must have $A_\D u_1\cap K[x_{r+1}, \ldots , x_n]\neq
0$. Pick a nonzero element, say $u_2$, of the intersection, then
applying several times maps of the type $\ad (\der_{x_k})$,
$r+1\leq k \leq n$, we have a nonzero scalar, hence $L=\CD (A_\D
)$. This proves that the algebra $\CD (A_\D )$ is a simple
algebra.
 $\Box $

\begin{corollary}\label{1CDAd=gd}%\marginpar{1CDAd=gd}
The derivations $\D(\i , \j )^{-1}\der_{\i ; \j , j_{r+1}}, \ldots
, \D(\i , \j )^{-1}\der_{\i ; \j , j_n}$ from Theorem
\ref{CDAd=gd} are respectively the partial derivatives
$\der_{j_{r+1}}:=\frac{\der}{\der x_{j_{r+1}}}, \ldots ,
\der_{j_n}:=\frac{\der}{\der x_n}$ of the algebra $A_\D $ (and of
the field of fractions of $A$).
\end{corollary}

{\it Remark}. The equality $\der_{j_{r+k}}:=\frac{\der}{\der
x_{j_{r+k}}}$ means that the derivation $\der_{j_{r+1}}$ is a
unique extension of the partial derivative $\frac{\der}{\der
x_{j_{r+1}}}$ of the polynomial algebra  $K[x_{r+1}, \ldots ,
x_n]$ to the algebra $A_\D$.

\begin{corollary}\label{2CDAd=gd}%\marginpar{2CDAd=gd}
 Let $Q$ be the field of fractions of the algebra $A$. Under the
 assumption of Theorem \ref{CDAd=gd}, the algebra $\CD (Q)$ is
 generated by the field $Q$ and the derivations
 $\der_{\i ; \j , j_{r+1}}, \ldots , \der_{\i ; \j ,
j_n}$ and $\Der_K(Q)=\bigoplus_{k=r+1}^nQ\der_{\i ; \j , j_k}$.
\end{corollary}

{\it Proof}. Note that $\CD (Q)\simeq Q\t_{A_\D}\CD(A_\D )$ and
$\Der_K(Q)\simeq Q\t_{A_\D }\Der_K(A_\D )$, and the results follow
from Theorem \ref{CDAd=gd}.  $\Box $

As a corollary we have a short (new) proof of the following key
result on differential operators on a regular algebra.

\begin{theorem}\label{3Feb05}%\marginpar{3Feb05}
Let the algebra $A$ be a regular algebra. Then the algebra $\CD
(A)$ of differential operators on $A$ is a simple Noetherian
algebra generated by $A$ and $\Der_K(A)$.
\end{theorem}

{\it Proof}. Let $\D$ be the subalgebra of ${\rm End}_K(A)$
generated by $A$ and $\Der_K(A)$. By Theorem \ref{CDAd=gd}, $\CD
(A)_{\D (\i , \j )}=\D_{\D (\i , \j )}$ for all non-singular $\i$
and $\j $, or equivalently $\prod_{\i , \j }\CD (A)_{\D (\i , \j
)}\t_{\CD (A)}(\CD (A)/\D ) =0$. By Theorem \ref{CDAffl}, we must
have $\CD (A)=\D $. By Theorems \ref{CDAffl} and \ref{CDAd=gd},
$\CD (A)$ is a simple Noetherian algebra.  $\Box $

\begin{proposition}\label{P3Feb05}%\marginpar{P3Feb05}
Let $\i , \i'\in \II_r$, $\j =(j_1, \ldots , j_r)\in \JJ_r$, $\j'
=(j_1', \ldots , j_{r+1}')\in \JJ_{r+1}$, and $\{ j_{r+1}, \ldots
j_n\} =\{ 1, \ldots , n\} \backslash \{ j_1, \ldots , j_r\}$. Then
\begin{enumerate}
 \item

\begin{equation}
 \der_{\i',\j'}=\D (\i , \j )^{-1}\sum_{l=1}^s(-1)^{r+1+\nu_l}\D (\i';
j_1', \ldots , \widehat{j_{\nu_l}'}, \ldots , j_{r+1}')\der_{\i
;\j , j_{\nu_l}'}
\end{equation}
where  $j_{\nu_1}', \ldots , j_{\nu_s}'$ are  the elements
 of the set $\{ j_1', \ldots , j_r', j_{r+1}'\} \backslash \{ j_1,
 \ldots , j_r\}$.
\item $\der_{\i',\j'}=(-1)^{r+1+k}\frac{\D (\i'; j_1', \ldots ,
\widehat{j_k'}, \ldots , j_{r+1}' )}{\D (\i ; j_1', \ldots ,
\widehat{j_k'}, \ldots , j_{r+1}'  )}\der_{\i ,\j'}$ provided $\D
(\i ; j_1', \ldots , \widehat{j_k'}, \ldots , j_{r+1}'  )\neq 0$.
\end{enumerate}
\end{proposition}

{\it Proof}. $1$. By Theorem \ref{25July04},  $\der_{\i' , \j'}=\D
(\i , \j )^{-1}\sum_{k=r+1}^n\l_k\der_{\i ; \j , j_k}$ for some
$\l_k\in A$. For each $k'=r+1, \ldots , n$, $\der_{\i ; \j ,
j_k}(x_{j_{k'}})=\d_{k,k'}\D (\i , \j )$. Then evaluating the
equality above at $x_{j_k}$, one
 gets the equality $\l_k= \der_{\i' ,\j'}(x_{j_k})$. So, $\l_k=0$
 if $j_k\not\in \{ j_{\nu_1}', \ldots , j_{\nu_s}'\}$, and if
 $j_k=j_{\nu_l}'$ then $\l_k=(-1)^{r+1+\nu_l}
 \D (\i ; j_1', \ldots , \widehat{j_{\nu_l}'}, \ldots ,
 j_{r+1}')$. This finishes the proof of the first statement.

 $2$. By the first statement where we put $\j =
 (j_1', \ldots , \widehat{j_k'}, \ldots ,
 j_{r+1}')$, we have $\der_{\i',\j'}=(-1)^{r+1+k}\frac{\D (\i'; j_1', \ldots ,
\widehat{j_k'}, \ldots , j_{r+1}' )}{\D (\i ; j_1', \ldots ,
\widehat{j_k'}, \ldots , j_{r+1}'  )}\der_{\i ,\j'}$.
 $\Box $

{\it Remark}. Let us fix elements $\i\in \II_r$ and $\j\in \JJ_r$.
Then, for each $\i'\in \II_r$ and $\j'\in \JJ_{r+1}$ (as above),
let $a(\i',\j')$ be the vector of coefficients $(\l_{r+1}, \ldots
, \l_n)$ from the proof of Proposition \ref{P3Feb05}. If $A$ is a
{\em regular} algebra, the vectors $a(\i',\j')$ form a generating
set for the $A$-module of solutions to the system of inclusions
from Theorem \ref{25July04} (by Theorem \ref{23July04}).

{\bf  Proof of Theorem \ref{9Feb05}}. Let $D=D(A)$ be the algebra
generated by $A$ and $d_{\i , \j }$ that satisfy the defining
relations (\ref{RD2}) and (\ref{RD3}).  By (\ref{RD2}), the inner
derivation $\ad (d_{\i , \j })$ of the algebra $D$ acts on the
algebra $A$ as the derivation $\der_{\i , \j }$. There is  a
natural $K$-algebra epimorphism $D\ra \CD =\CD (A)$ which maps
$a\mapsto a$ $(a\in A)$ and $d_{\i , \j }\mapsto \der_{\i \, \j
}$. Since $A$ is a subalgebra of $\CD $, the algebra $A$ is a
subalgebra of $D$ as well. It follows from (\ref{RD2}) that
$D=\sum A d_*\cdots d_*$ where $d_*$ stands for any element $d_{\i
, \j }$, and then  for any element $0\neq s\in A$   there exists a
left and right Ore localization of the algebra $D$ at the powers
of the element $s$, denoted $D_s$, and $D_s=\sum A_sd_*\cdots
d_*$.
 We have a commutative diagram of algebra homomorphisms

$$\xymatrix{D\ar[r]\ar[d]&\prod_{\i\in \II_r , \j \in \JJ_r}D_{\D (\i , \j )}\ar[d]\\
\CD \ar[r]&\prod_{\i \in \II_r, \j\in \JJ_r }\CD_{\D (\i , \j )}}
$$
where horizontal maps are faithfully flat extensions (since $A=(\D
(\i, \j ))_{\i\in \II_r, \j \in \JJ_r}$ as $A$ is regular) and
vertical maps are epimorphisms. By Theorem \ref{CDAd=gd}.(1), Proposition \ref{P3Feb05}, (\ref{1Derel}) and
(\ref{RD3}), the epimorphism $D_{\D (\i , \j )}\ra \CD_{\D (\i ,
\j )}$ is in fact an isomorphism. So, the right vertical map is an
algebra isomorphism. Then the
 left vertical map must be an isomorphism (by faithfully flatness), i.e. $D\simeq \CD $.
 $\Box $

{\bf  Proof of Theorem \ref{9bFeb05}}. The arguments are similar
to that in the proof of Theorem \ref{9Feb05}. Let ${\rm DER}(A)$
be a left  $A$-module generated by symbols $\der_{\i , \j }$
subject to the defining relations (\ref{Derel}). We have a
commutative diagram of left $A$-modules:
$$\xymatrix{{\rm DER}(A)\ar[r]\ar[d]&\prod_{\i\in \II_r , \j \in \JJ_r}{\rm DER}(A)_{\D (\i , \j )}\ar[d]\\
\Der_K(A) \ar[r]&\prod_{\i \in \II_r, \j\in \JJ_r }\Der_K(A)_{\D
(\i , \j )}}
$$
where horizontal maps are faithfully flat $A$-module monomorphisms
 as  $A=(\D (\i , \j ))_{\i\in \II_r, \j \in \JJ_r}$
 ($A$ is regular) and vertical maps are natural $A$-module
 epimorphisms. By Proposition  \ref{P3Feb05}, Theorem \ref{25July04},  and (\ref{Derel}),
 each epimorphism ${\rm DER}(A)_{\D (\i , \j )}\ra \Der_K(A)_{\D (\i ,
 \j )}$ is an isomorphism. So, the right vertical map must be an
 isomorphism, and so the left vertical map must be an
 isomorphism (by faithfully flatness), i.e. ${\rm DER}(A)\simeq \Der_K(A)$.   $\Box $

{\bf Proof of Theorem \ref{25Feb05}}. $(1\Rightarrow 2)$ By
Theorem \ref{9Feb05}, $\CD (A) = \gD (A)$ and the result follows
from Theorem \ref{MR15.3.8} (or from Theorem \ref{3Feb05}).

$(2\Rightarrow 3 )$ Suppose that the $\gD (A)$-module $A$ is not
simple then it contains a proper ideal, say $\ga$, stable under
${\rm der}_K(A)$. Then $\ga \gD (A)$ is a proper  ideal of the
algebra $\gD (A)$ since $0\neq \ga \gD (A)(A)\subseteq \ga$, a
contradiction.

$(3\Rightarrow 1)$ By Lemma \ref{sid}, ${\rm der}_K(A)(A)=\ga_r$
is a nonzero $\gD (A)$-submodule of $A$, therefore $A=\ga_r$ since
$A$ is a simple $\gD (A)$-module. So, $A$ is a regular algebra.
$\Box $

Any element $\d $ of the algebra $\gD (A)$ is a sum of elements
$ad$ where $a\in A$ and $d$ is a product of the derivations
$\derij$.  So, the algebra $\gD(A)=\cup_{i\geq 0}\gD(A)_i$ has a
natural filtration by the total degree of natural derivations, and
the associated graded algebra ${\rm gr} \, \gD
(A):=\bigoplus_{i\geq 0}\gD (A)_i/\gD (A)_{i-1}$ is a {\em
commutative finitely generated} algebra. Therefore, {\em the
algebra $\gD (A)$ is a (left and right) Noetherian algebra}. The
algebra $\gD (A)$ is a {\em domain} as a subalgebra of the domain
$\CD (Q)$. $\Kdim$ stands for the  (left or right) {\em Krull
dimension} (they both coincide for algebras we consider). For  a
definition and properties of $\Kdim$ the reader is referred to
\cite{MR}, Ch. 6.

%*** Qu: $\gD (A)_i=\gD (A)\cap \CD (A)_i$?   ***
%*** Qu: Is $\gr (\gD (A))$ a domain?  ***
%***  characteristic var. of gr $\CD (A)$, sing of $A$ coincide
%with sing of $gr \CD (A)$  ***
\begin{lemma}\label{KdA}%\marginpar{KdA}
\begin{enumerate}
\item The algebra $\gD (A)$ is a (left and right) Noetherian
domain.  \item $\GK (\gD (A)) =2\Kdim (A)$. \item $\Kdim (A)\leq
\Kdim (\gD (A))\leq 2\Kdim (A)$.\item If $A$ is a regular algebra
then $\Kdim (\gD (A))={\rm gldim} (\gD (A))= \frac{1}{2}\GK (\gD
(A))=\Kdim (A)$.
\end{enumerate}
\end{lemma}

{\it Proof}. $2$. Let $d=\Kdim (A)$. Recall that $\GK (\D (A))=2d$
(\cite{MR}, 15.3.6.(ii)). Then $\gD (A)\subseteq \D (A)$ implies
$\GK (\gD (A))\leq \GK (\D (A))=2d$.

Fix $\i \in \II_r$ and $\j\in\JJ_r$ then $\D =\D(\i, \j )\neq 0$
and the localization $A_\D$ of the algebra $A$ at the powers of
the element $\D $ is a regular algebra, and so the algebra  $\gD
(A_\D )=\CD (A_\D )\simeq \CD (A)_\D $ (Theorem \ref{9Feb05})
contains the Weyl algebra $A_{n-r}=K\langle x_{j_{r+1}}, \ldots ,
x_{j_n}, \der_{j_{r+1}}, \ldots , \der_{j_n}\rangle $ (see
Corollary \ref{1CDAd=gd} and Theorem \ref{CDAd=gd}). It follows
immediately that standard monomials in $ x_{j_{r+1}}, \ldots ,
x_{j_n}, \der_{\i; \j , j_{r+1}}:= \D\der_{j_{r+1}}, \ldots ,
\der_{\i; \j , j_n}:= \D\der_{j_n}\in \gD (A)$ are $K$-linearly
independent in $\gD (A)$ (a standard monomial means an element of
the type $x_{j_{r+1}}^{n_{r+1}} \cdots x_{j_n}^{n_n}\der_{\i; \j ,
j_{r+1}}^{m_{r+1}} \cdots  \der_{\i; \j , j_n}^{m_n}$ where
$n_{r+1} , \ldots , m_n$ are non-negative integers). Therefore,
$\GK (\gD (A))\geq \GK (A_{n-r})=2(n-r)=2d$, and so $\GK (\gD
(A))=2d$.

$3$. $\gD (A)_\D\simeq \gD (A_\D )=\CD (A_\D )$ implies $\Kdim
(\gD (A))\geq \Kdim (\CD(A_\D ))=\Kdim (A_\D )=d$.  Recall that
the algebra $\gD (A)=\cup_{i\geq 0}\gD (A)_i$ is the filtered
algebra such that the associated graded algebra ${\rm gr} (\gD
(A))$ is a commutative finitely generated algebra, and so $\Kdim (
{\rm gr} \, \gD (A))=\GK ( {\rm gr} \, \gD (A))$ (\cite{MR},
8.2.14).
\begin{eqnarray*}
\Kdim (\gD (A))&\leq & \Kdim ({\rm gr} \, \gD (A)) \;\;\;
(\cite{MR}, 6.5. 6)\\
&=& \GK ({\rm gr} \, \gD (A))\leq \GK (\gD (A))=2d \;\;
(\cite{MR}, 8.3.20).
\end{eqnarray*}

$4$. If $A$ is regular then $\gD (A)=\CD (A)$, and for the algebra
$\CD (A)$ the statements are well-known (see \cite{MR}, Chapter
15). $\Box $

%%%%%%%%%%%%%%%%%% SECTION  %%%%%%%%%%%%%%%%%%%%%%%%

\section{Generators and defining relations for the ring of differential operators
on a regular algebra of essentially finite
type}\label{essgrel}%\marginpar{essgrel}

In this section, we prove that the main results for the algebra
$A$ (Theorems \ref{9bFeb05}, \ref{9Feb05}, and \ref{25Feb05}) also
hold for localizations of the algebra $A$.

{\it Definition}. A localization of an {\em affine} algebra is
called an algebra of {\bf essentially finite type}.

Let $\CA :=\S1 A$ be a localization of the algebra $A=P_n/I$ at a
multiplicatively closed subset $S$ of $A$.  Recall that the
functors $\Der_K(\cdot )$ and $\CD (\cdot )$ commute with
localizations: $\Der_K(\S1 A)\simeq \S1 \Der_K(A)$, $\CD (\S1
A)\simeq \S1 \CD (A)$ and $\CD (\S1 A)_i\simeq \S1 \CD (A)_i$ for
all $i\geq 0$.

{\it Definition}. The $\CA$-module ${\rm der}_K(\CA ):= \sum_{\i
\in \II_r, \j\in \JJ_{r+1}}\CA \derij$  is called the $\CA$-module
of {\bf natural derivations} of the algebra $\CA$ of essentially
finite type, and the subalgebra of $\CD (\CA )$ generated by $\CA
$ and ${\rm der}_K(\CA )$ is called the algebra of {\bf natural
differential operators} on $\CA$ denoted  $\gD (\CA )$.

By the very definition,   the functors ${\rm der}_K(\cdot )$ and
$\gD (\cdot )$ commute with localizations: ${\rm der}_K(\S1
A)\simeq \S1 {\rm der}_K(A)$, $\gD (\S1 A)\simeq \S1 \gD (A)$ and
$\gD (\S1 A)_i\simeq \S1 \gD (A)_i$ for all $i\geq 0$.

The next lemma follows from these facts and from Theorem
\ref{7Mar05}.

\begin{lemma}\label{ess7Mar05}%\marginpar{ess7Mar05}
\begin{enumerate}
\item The $\CA$-module ${\rm der}_K(\CA )$ does not depend on the
choice of presentation of the algebra $\CA$ as $\S1 P_n/I$ (i.e.
$\S1 P_n/I\simeq S'^{-1} P_{n'}/I'$ implies ${\rm der}_K(\S1
P_n/I)\simeq {\rm der}_K(S'^{-1} P_{n'}/I')$). \item The algebra
of natural differential operators $\gD (\CA )$ on $\CA$   does not
depend on the presentation of the algebra $\CA$ as $\S1 P_n/I$.
\end{enumerate}
\end{lemma}

The {\em Jacobian criterion of regularity} holds for the algebra
 $\CA$: {\em the algebra $\CA$ is a regular algebra iff} $\S1
\ga_r=\CA$ (\cite{Ma}, Theorem 30.3 and Remark 2, p.235). So, if
the algebra $\CA$ is {\em regular} we have the faithfully flat
extension of left $\CA$-modules:
$$\Der_K(\CA )\ra \prod_{\i\in \II_r , \j \in \JJ_r}\Der_K(\CA )_{\D (\i , \j
)},$$ and the faithfully flat extension of algebras
$$\CD (\CA )\ra \prod_{\i\in \II_r , \j \in \JJ_r}\CD (\CA )_{\D (\i , \j
)}.$$ These are the main ingredients of the proofs of Theorems
\ref{9Feb05} and \ref{9bFeb05}. Now, the analogous results for the
algebra $\CD (\CA )$ follow easily using the facts that that
differential operators and natural differential operators commute
with localizations.

\begin{theorem}\label{ess9bFeb05}%\marginpar{ess9bFeb05}
Let the algebra $\CA$ be a regular algebra. Then the left
$\CA$-module $\Der_K(\CA )$ is generated by derivations $\der_{\i
, \j }$, $\i \in \II_r$, $\j \in \JJ_{r+1}$ that satisfy the
defining relations (\ref{Derel}) (as a left $\CA $-module).
\end{theorem}

{\it Proof}. Repeat the proof of Theorem \ref{9bFeb05}.  $\Box$

\begin{theorem}\label{ess9Feb05}%\marginpar{ess9Feb05}
Let the algebra $A$ be a regular algebra. Then the ring of
differential operators $\CD (\CA )$ on $\CA$ is generated over $K$
by the algebra $\CA$ and the derivations  $\der_{\i , \j }$, $\i
\in \II_r$, $\j \in \JJ_{r+1}$ that satisfy the defining relations
(\ref{Derel}) and (\ref{1Derel}). 
\end{theorem}

{\it Proof}. It is well-known that the ring of differential
operators on a regular domain of essentially finite type is
generated by the algebra itself and its derivations. By Theorem
\ref{ess9bFeb05}, the algebra $\CD (\CA )$ is generated by $\CA$
and $\der_{\i , \j }$, $\i \in \II_r$, $\j \in \JJ_{r+1}$. It is
well-known that the algebra $\CD (\CA )$ is a simple Noetherian
domain. Now, the facts just mentioned  and proved and using the
faithfully flat extension of algebras $\CD (\CA )\ra \prod_{\i\in
\II_r , \j \in \JJ_r}\CD (\CA )_{\D (\i , \j )}$, one can repeat
word for word the proof of Theorem \ref{9Feb05} to finish the
proof of theorem.  $\Box $.

\begin{theorem}\label{ess25Feb05}%\marginpar{ess25Feb05}
(Criterion of regularity via $\gD (\CA )$) The following
statements are equivalent.
\begin{enumerate}
\item $\CA$ is a regular algebra. \item $\gD (\CA )$ is a simple
algebra.\item $\CA$ is a simple $\gD (\CA )$-module.
\end{enumerate}
\end{theorem}

{\it Proof}. $(1\Rightarrow 2)$ If $\CA$ is a regular algebra then
$\CD (\CA )$ is a simple Noetherian algebra generated by $\CA$ and
$\Der_K(\CA )$. By Theorem \ref{ess9Feb05}, $\Der_K(\CA )={\rm
der}_K(\CA )$, hence $\gD (\CA ) = \CD (\CA )$ is a simple
algebra.

$(2\Rightarrow 3 )$ Suppose that the $\gD (\CA )$-module $\CA$ is
not simple then it contains a proper ideal, say $\ga$, stable
under ${\rm der}_K(\CA )$. Then $\ga \gD (\CA )$ is a proper ideal
of the algebra $\gD (\CA )$ since $0\neq \ga \gD (\CA )(\CA
)\subseteq \ga$, a contradiction.

$(3\Rightarrow 1)$ Lemma \ref{sid} holds for the algebra $\CA$ as
well (with the same proof), and so ${\rm der}_K(\CA )(\CA )=\S1
\ga_r$ is a nonzero $\gD (\CA )$-submodule of $\CA$, therefore
$\CA =\S1 \ga_r$ since $\CA$ is a simple $\gD (\CA )$-module. So,
$\CA$ is a regular algebra. $\Box $

 The algebra $\gD(\CA )=\cup_{i\geq 0}\gD(\CA )_i$ has a
natural filtration by the total degree of natural derivations, and
the associated graded algebra ${\rm gr} \, \gD (\CA
):=\bigoplus_{i\geq 0}\gD (\CA )_i/\gD (\CA )_{i-1}$ is a {\em
commutative finitely generated $\CA$-algebra}.
\begin{lemma}\label{essKdA}%\marginpar{essKdA}
\begin{enumerate}
\item The algebra $\gD (\CA )$ is a (left and right) Noetherian
domain.  \item $\GK (\gD (\CA )) =2\GK (\CA )$. \item $\Kdim (\CA
)\leq \Kdim (\gD (\CA ))\leq 2\GK (\CA)$.
\end{enumerate}
\end{lemma}

{\it Proof}. $1$. Since $\gD (\CA )=\gD (\S1 A)\simeq \S1 \gD
(A)$, the algebra $\gD (\CA )$ is a Noetherian domain (Lemma
\ref{KdA}. (1)).

$2$. $A\subseteq \CA \subseteq Q$ implies (by Lemma \ref{KdA}.(2))
\begin{eqnarray*}
 2\GK (\CA )&=&2\GK (A)=2\Kdim (A)=\GK (\gD (A))\leq \GK (\gD (\CA
))\\
&\leq &\GK (\CD (Q))=2\GK (Q) =2\GK (\CA ),
\end{eqnarray*}
and so $\GK (\gD (\CA ))=2\GK (\CA )$.

$3$.  Fix $\i \in \II_r$ and $\j\in\JJ_r$ then $\D =\D(\i, \j
)\neq 0$, and the localization $\CA_\D$ of the algebra $\CA$ at
the powers of the element $\D $ is a regular algebra. $ \gD (\CA
)_\D \simeq \gD (\CA_\D )=\CD (\CA_\D )$ (Theorem \ref{ess9Feb05})
implies $\Kdim (\gD (\CA ))\geq \Kdim (\CD (\CA_\D ))\geq \Kdim
(\CA_\D )$, the last inequality is due to the fact that the map
$U\mapsto \CD (\CA_\D )U$ (resp. $U\mapsto U\CD (\CA_\D )$) from
the set of ideals of the algebra $\CA_\D$ to the set of left
(resp. right) ideals of the algebra $\CD (\CA_\D )$ is an
injection since $\CD (\CA_\D )=\CA_\D \oplus \CD (\CA_\D )\Der_K(
\CA_\D )$ (resp. $\CD (\CA_\D )=\CA_\D \oplus \Der_K( \CA_\D )\CD
(\CA_\D )$). Now,
$$ \Kdim (\CD (\CA ))=\Kdim( \S1 \CD (A))\leq \Kdim (\CD (A))\leq
2\Kdim (A)=2\GK (A)=2\GK (\CA ),$$ by Lemma \ref{KdA}.(3). $\Box $

{\it Remark}. A typical situation when the results of this section
 can be applied is the ring of differential operators on the algebra
of regular functions in the neighbourhood of a {\em simple} point
(of a singular variety): given a prime ideal $\gp$ of $A$ such
that $\ga_r\not\subseteq \gp$, then $A_\gp$ is a smooth domain of
essentially finite type, and so for the ring of differential
operators $\CD (A_\gp )$ all the results of this section hold.

%%%%%%%%%%%%%%%%%% SECTION  %%%%%%%%%%%%%%%%%%%%%%%%

\section{Ring of differential operators
 on a singular irreducible affine algebraic
variety}\label{opsing}%\marginpar{opsing}

 \begin{lemma}\label{duchi}%\marginpar{duchi}
Let $R=K\langle x_1, \ldots , x_n\rangle$ be a commutative
finitely generated algebra over the field $K$, and $\CD (R)$ be
the ring of differential operators on $R$. Each element $\d \in
\CD (R)_i$ is completely determined by its values on the elements
$x^\alpha$, $\alpha \in \mathbb{N}^n$, $| \alpha |\leq i$.
\end{lemma}

{\it Proof}. It suffices to prove that given $\d \in \CD (R)_i$
satisfying $\d (x^\alpha )=0$ for all $\alpha$ such that $| \alpha
| \leq i$ then $\d =0$. We use induction on $i$. The case $i=0$ is
trivial: $\d \in \CD (R)_0=R$ and $0=\d \cdot 1=\d$. Suppose that
the statement is true for all $i'<i$. For each $x_j$, $[\d ,
x_j]\in \CD (R)_{i-1}$ and, for each $x^\alpha$ with $|\alpha |
\leq i-1$, $[\d , x_j](x^\alpha )=\d (x_jx^\alpha)-x_j\d
(x^\alpha)=0$. By induction, $[\d , x_j]=0$. Now, for any
$x^\alpha$,
$$ \d (x^\alpha )=\d (x_jx^{\alpha -e_j})=x_j\d (x^{\alpha
-e_j})+[\d , x_j](x^{\alpha -e_j})=x_j\d (x^{\alpha -e_j})=\cdots
= x^\alpha \d (1)=0,$$
 and so $\d =0$, as required.
$\Box $

Let $S$ be a multiplicatively closed subset of the algebra $A$.
Let us consider the natural inclusion  $\CD (A)\subseteq \S1 \CD
(A)$ of filtered algebras (by the total degree of derivations).
$\CD (A)=\{ \d \in \S1 \CD (A)\, | \, \d (A)\subseteq A\}$ and
$\CD (A)_i=\{ \d \in \S1 \CD (A)_i\, | \, \d (A)\subseteq A\}$,
$i\geq 0$.

\begin{lemma}\label{1duchi}%\marginpar{1duchi}
Given $\d \in \S1 \CD (A)_i$. Then $\d \in \CD (A)$ iff $\d
(x^\alpha )\in A$ for all $\alpha$ such that $| \alpha |\leq i$.
\end{lemma}

{\it Proof}. $(\Rightarrow )$ Trivial.

$(\Leftarrow )$  We use induction on $i$. When  $i=0$, $\d \in \S1
\CD (A)_0=\S1 A$ and $\d =\d (1)=\d \cdot 1\in A$.  Suppose that
the statement is true for all $i'<i$. Let $\d \in \S1 \CD (A)_i$
satisfy $\d (x^\alpha )\in A$ for all $\alpha $ such that $|
\alpha | \leq i$. For each $j$,
 $[\d ,
x_j]\in \S1 \CD (R)_{i-1}$ and, for each $x^\alpha$ with $|\alpha
| \leq i-1$, $[\d , x_j](x^\alpha )=\d (x_jx^\alpha)-x_j\d
(x^\alpha)\in A$, and so, by induction, $[\d , x_j]\in \CD
(A)_{i-1}$. Now, for any $x^\alpha$,
\begin{eqnarray*}
\d (x^\alpha )&=&\d (x_jx^{\alpha -e_j})=x_j\d (x^{\alpha
-e_j})+[\d , x_j](x^{\alpha -e_j})\equiv x_j\d (x^{\alpha -e_j})\,
{\rm mod} \, A\\
&\equiv &\cdots \equiv  x^\alpha \d (1)\equiv 0 \, {\rm mod} \, A,
\end{eqnarray*}
 and so $\d \in \CD (A)_i$.
$\Box $

\begin{proposition}\label{sinDAfg}%\marginpar{sinDAfg}
Suppose that $\i =(i_1, \ldots , i_r)$ and $\j =(j_1, \ldots ,
j_r)$ are non-singular, let $\der_{j_{r+1}}:=\D^{-1}\der_{\i ; \j,
j_{r+1}}, \ldots ,\der_{j_n}:=\D^{-1}\der_{\i ; \j, j_n}$ where
$\D =\D (\i , \j )$ and $\{ j_{r+1}, \ldots , j_n\} =\{ 1, \ldots
, n\} \backslash \{ j_1, \ldots , j_r\}$.
\begin{enumerate}
\item For each $i\geq 0$, $\CD (A)_i=\{ \d \in \sum_{| \alpha
|\leq i} A\der^\alpha\, | \,  \d (x^\beta )\in A$ for all $\beta
\in \mathbb{N}^n\}$ where $\alpha =(\alpha_{r+1}, \ldots ,
\alpha_n)\in \mathbb{N}^{n-r}$ and  $\der^\alpha
=\prod_{k=r+1}^n\der_{j_k}^{\alpha_k}$ (recall that the
derivations $\der_{j_{r+1}}, \ldots , \der_{j_n}$ commute). \item
For each $i\geq 0$, $\CD (A)_i$ is a finitely generated left
$A$-module.
\end{enumerate}
\end{proposition}

{\it Proof}. $1$. Clearly, $\CD (A)_i=\{ \d \in \sum_{|\alpha
|\leq i}A_\D \der^\alpha\, | \, \d (A)\subseteq A\}$. Then, by
Lemma \ref{1duchi}, $\CD (A)_i=\{ \d \in \sum_{|\alpha |\leq
i}A_\D \der^\alpha\, | \, \d (x^\beta)\subseteq A$ for all $\beta
\in \mathbb{N}^n$ such that $|\beta | \leq i\}$. Note that
$\der_{j_{r+1}}=\frac{\der}{\der x_{j_{r+1}}}, \ldots ,
\der_{j_n}=\frac{\der}{\der x_{j_n}}$. Then the conditions that
$\d (x_{j_{r+1}}^{\gamma_{r+1}}\cdots x_{j_n}^{\gamma_n})\in A$
for all $\gamma=(\gamma_{r+1}, \ldots ,\gamma_n)\in
\mathbb{N}^{n-r}$ with $|\gamma | \leq i$ are equivalent to $\d
\in \sum_{|\alpha | \leq i}A\der^\alpha$. This gives the first
statement.

$2$. $\CD (A)_i$ is a Noetherian left $A$-module as a submodule of
the Noetherian $A$-module $\sum_{|\alpha | \leq i}A\der^\alpha$,
and so  $\CD (A)_i$ is a finitely generated left $A$-module.
$\Box $

Department of Pure Mathematics

University of Sheffield

Hicks Building

Sheffield S3 7RH

UK

email: v.bavula@sheffield.ac.uk

\end{document}